\documentclass{amsart}
\usepackage{graphicx}
\usepackage{amssymb}
\usepackage{epstopdf}
\usepackage{nicefrac}
\usepackage{enumerate}
\usepackage{float}

\usepackage{tabularx}
\usepackage{hyperref}

\newtheorem{thm}{THEOREM}[section]

\newtheorem{cor}[thm]{COROLLARY}
\newtheorem{defn}[thm]{DEFINITION}

\newtheorem{lemma}[thm]{LEMMA}
\newtheorem{prob}[thm]{PROBLEM}
\newtheorem{prop}[thm]{PROPOSITION}

\input xy
\xyoption {all}

\newcommand{\ds}{\displaystyle}
\newcommand{\vp}{{\varphi}}

\newcommand{\bK}{{\bf K}}
\newcommand{\cE}{{\mathcal E}}
\newcommand{\cG}{{\mathcal G}}
\newcommand{\cGF}{\cG_{\F}} 
\newcommand{\cGX}{\cG_{\fX}} %
\newcommand{\cGY}{\cG_{\fY}} %
\newcommand{\cH}{{\mathcal H}}
\newcommand{\cI}{{\mathcal I}}
\newcommand{\cL}{{\mathcal L}}
\newcommand{\cP}{{\mathcal P}}
\newcommand{\cQ}{{\mathcal Q}}
\newcommand{\cR}{{\mathcal R}}
\newcommand{\cRF}{\cR_{\F}} 
\newcommand{\cS}{{\mathcal S}}
\newcommand{\cT}{{\mathcal T}}
\newcommand{\cU}{{\mathcal U}}
\newcommand{\cV}{{\mathcal V}}
\newcommand{\cX}{{\mathcal X}}
\newcommand{\cZ}{{\mathcal Z}}
\newcommand{\dA}{d_{\fA}} 
\newcommand{\dF}{d_{\F}} 
\newcommand{\dFU}{\delta^{\F}_{\cU}} 
\newcommand{\dM}{d_{\fM}} 
\newcommand{\Dom}{{\rm Dom}} %
\newcommand{\dX}{d_{\fX}} 
\newcommand{\dXp}{d_{\fX}'} 
\newcommand{\dY}{d_{\fY}} 
\newcommand{\e}{{\epsilon}} 
\newcommand{\eFU}{{\epsilon^{\F}_{\cU}}} 
\newcommand{\eM}{{\epsilon_{\fM}}} 
\newcommand{\eU}{{\epsilon_{\cU}}} 
\newcommand{\F}{{\mathcal F}} 
\newcommand{\fA}{{\mathfrak{A}}}
\newcommand{\fM}{{\mathfrak{M}}}
\newcommand{\fN}{{\mathfrak{N}}}

\newcommand{\fX}{{\mathfrak{X}}}
\newcommand{\fY}{{\mathfrak{Y}}}
\newcommand{\G}{\Gamma}
\newcommand{\GF}{\Gamma_{\F}} 
\newcommand{\lF}{{\lambda_{\mathcal F}}} 
\newcommand{\mF}{{\mathbb F}}
\newcommand{\mG}{{\mathbb G}}
\newcommand{\mH}{{\mathbb H}}
\newcommand{\mR}{{\mathbb R}}
\newcommand{\mS}{{\mathbb S}}
\newcommand{\mT}{{\mathbb T}}
\newcommand{\mZ}{{\mathbb Z}}
\newcommand{\oU}{{\overline{U}}}
\newcommand{\ovx}{\overline{x}} %
\newcommand{\ovy}{\overline{y}} %
\newcommand{\whU}{{\widehat U}}
\newcommand{\whvarp}{{\widehat \varphi}}
\newcommand{\wtM}{\widetilde{M}}

\begin{document}

 \begin{abstract}
A matchbox manifold is a connected, compact foliated space with totally disconnected transversals; or in other notation, a generalized lamination. It is said to be Lipschitz if there exists a metric on its   transversals for which the holonomy maps are Lipschitz. Examples of Lipschitz matchbox manifolds include the exceptional minimal sets for $C^1$-foliations of compact manifolds, tiling spaces, the classical solenoids, and the weak solenoids of McCord and Schori, among others.  We address the question: When does a Lipschitz matchbox manifold admit an embedding as a minimal set for a smooth dynamical system, or more generally for as an exceptional minimal set for a $C^1$-foliation of a smooth  manifold? We gives examples which do embed, and develop criteria for showing when they  do not embed, and give examples. We also discuss the classification theory for Lipschitz weak solenoids.
\end{abstract}

\title{Lipschitz  matchbox manifolds}


\author{Steven Hurder}
\address{Steven Hurder, Department of Mathematics, University of Illinois at Chicago, 322 SEO (m/c 249), 851 S. Morgan Street, Chicago, IL 60607-7045}
\email{hurder@uic.edu}

\thanks{Version date: September 5, 2013; updated August 26, 2014}

\date{}


\maketitle


\section{Introduction} \label{sec-intro}

 In this paper, we consider the classification and embedding problems for matchbox manifolds, from the viewpoint of Lipschitz pseudogroups, and develop  invariants which are obstructions to realizing a matchbox manifold as a minimal set. Matchbox manifolds are a class of continua that occur naturally in the study of  dynamical systems, and in foliation theory as   exceptional minimal sets.  The overall goal of this research program is to develop tools for the classification of these spaces, and to understand which matchbox manifolds are homeomorphic to exceptional minimal sets for  $C^r$-foliations, for $r \geq 1$.

We first discuss an important  motivation for   interest in this program of study. 
\begin{prob}[Sondow \cite{Sondow1975}] \label{problem4}
 When is a smooth connected $n$-manifold $L$ without boundary, diffeomorphic to a leaf of a foliation $\F_M$ of a compact smooth manifold $M$? 
 \end{prob}
For the case where $L$ has dimension $n=1$,  the problem is trivial. Also, for dimension $n=2$,  Cantwell  and Conlon showed in \cite{CC1987}  that any surface without boundary is a diffeomorphic to a leaf of a smooth codimension-$1$  foliation of a compact $3$-manifold. 

On the other hand, Ghys   \cite{Ghys1985}  and   Inaba,  Nishimori, Takamura and  Tsuchiya   \cite{INTT1985}, constructed   $3$-manifolds  which are not homeomorphic to a  leaf of any codimension-$1$ foliation of a compact manifold.  Souza and Schweitzer \cite{SouzaSchweitzer2013} give further examples in higher dimensions, of manifolds which cannot be leaves in codimension one.  The non-embedding examples by these authors are essentially the only known results on   Problem~\ref{problem4} in this generality, and they are  for codimension-one foliations.  

There is  a natural variant of Problem~\ref{problem4},  posed in the 1974 ICM address by   Sullivan  \cite{Sullivan1975}:
\begin{prob}\label{problem5}
 Let $L$ be a complete Riemannian smooth manifold without boundary. When is   $L$ \emph{quasi-isometric} to a leaf of a $C^r$-foliation $\F_M$ of a compact smooth manifold $M$, for $r\geq 1$?
  \end{prob}
A quasi-isometric embedding of $L$ must preserve its quasi-invariant geometric properties, which can be used to construct obstructions to such an embedding. For example, Cantwell and Conlon studied in 
\cite{CC1977}, \cite{CC1978} how the asymptotic behavior of the metric on $L$ is related to the dynamics of the leaf in a codimension-one foliation.  

The work of Phillips and Sullivan in  \cite{PS1981} introduced the asymptotic Euler class of a non-compact Riemannian $2$-manifold $L$ which has  subexponential volume growth rate,   and showed this can be used as an obstruction to a quasi-isometric embedding of   $L$ as a leaf, depending on the topology of the ambient manifold $M$. This result was generalized by  Januszkiewicz in   \cite{Januszkiewicz1984} to obtain obstructions in terms of  the \emph{asymptotic Pontrjagin numbers} of an open Riemannian $n$-manifold with subexponential volume growth rate, for $n =4k$ with $k \geq 1$.

In an alternate direction, Attie and Hurder in \cite{AttieHurder1996} introduced an invariant of open manifolds, its ``leaf entropy'', or ``asymptotic leaf complexity'', and gave examples of open manifolds with exponential volume growth rate that cannot be quasi-isometric to a leaf in a foliation of any codimension.  Examples of surfaces with exponential growth rate that cannot be quasi-isometric to a leaf  were constructed by Schweitzer in \cite{Schweitzer1995} and Zeghib in \cite{Zeghib1994}, using a variant of the approach in \cite{AttieHurder1996}. 
The work of  Schweitzer      \cite{Schweitzer2011} exhibits further  examples of complete Riemannian manifolds which are not quasi-isometric  to a leaf in any codimension-one foliation. The work of the author and Lukina \cite{HL2014} generalizes the results of \cite{AttieHurder1996} to the broader class of matchbox manifolds.

The   non-embedding results mentioned above    rely on the simple strategy,  that a leaf in a compact foliated manifold $M$ has some type of recurrence properties,  and the idea is to formulate such a property, \emph{intrinsic} to $L$, which cannot be satisfied  if   $L$ is homeomorphic to a leaf, or possibly quasi-isometric to a leaf. 
Each such criteria for \emph{non-recurrence} then yields  non-embeddability results.

A leaf $L$  contained in a minimal set $\fM$ for a foliation $\F_M$ on a compact manifold $M$ has much stronger recurrence properties. 
For example,    Cass    observed in  \cite{Cass1985} that such a leaf must be ``quasi-homogeneous'', and that this property is an invariant of the quasi-isometry class of a Riemannian metric on $L$. He consequently gave examples of complete Riemannian manifolds, including leaves of foliations, which cannot be quasi-isometric to a leaf in a minimal set.  For example, Cass showed that any non-compact leaf in a Reeb foliation of $\mS^3$ cannot be realized as a leaf of a minimal set in any codimension.

The question raised by Cass' work suggests a   variant of the above questions, where we consider the closure $\fM = \overline{L}$ of a non-compact leaf  $L \subset M$, where   $\fM$    has the structure of a \emph{foliated space}. 
The formal definition  of a foliated space $\fM$   was given by Moore and Schochet \cite[Chapter 2]{MS2006}, as part of their development of a general formulation of the Connes measured leafwise-index theorem \cite{Connes1994}.    Candel and Conlon   \cite[Chapter 11]{CandelConlon2000} further developed the theory of foliated spaces, and gave many interesting examples. We are particularly interested in those cases where the transverse model space for the foliated space $\fM$ is  totally disconnected. 

A compact connected foliated space $\fM$ with totally disconnected transversals is  called a ``matchbox manifold'',  in accordance with   terminology introduced in continua theory \cite{AO1991,AO1995,AM1988}.   A matchbox manifold with $2$-dimensional leaves is a lamination by surfaces, as defined in \cite{Ghys1999,LM1997}.  If all leaves of $\fM$ are dense, then it is called a \emph{minimal matchbox manifold}.  A compact minimal set $\fM \subset M$ for a foliation $\F_M$ on a   manifold $M$ yields a foliated space with foliation $\F = \F | \fM$. If the minimal set is exceptional, then $\fM$ is   a minimal matchbox manifold. The formal definition and some basic   properties of matchbox manifolds are discussed in Section~\ref{sec-foliated}.

  The leaves of the foliation $\F$ of a foliated space $\fM$    admit a smooth Riemannian metric, and    for each leaf $L \subset \fM$ there is a well-defined quasi-isometry class of Riemannian metrics on $L$. 
      The   obstructions used in the works above, to show that a particular Riemannian manifold $L$ cannot be  quasi-isometric to a leaf of a foliation of a compact manifold $M$, also provide obstructions to realizing $L$ as a leaf in a compact foliated space $\fM$.

 The following problem is addressed in this work: 
 
  \begin{prob}\label{problem7}
 Let $\fM$ be a  minimal matchbox manifold. Does there exists a   homeomorphism  of $\fM$ to an exceptional  minimal  set of  a $C^r$-foliation $\F_M$ of a manifold $M$, for $r \geq 1$?
 \end{prob}

When such an embedding exists, then each leaf    $L \subset \fM$ is  quasi-isometric to a leaf of $\F_M$. If the leaf $L$ is dense in $\fM$ and $\fM$ is non-embeddable, then this gives   a criteria for the non-embedding of $L$, that depends not just on the intrinsic geometry and topology of  $L$, but includes   ``extrinsic properties'' of $L$ in $\fM$, such as  the   transverse geometry and dynamics of  the foliated space $\fM$. 


Observe that  if $\fM$ is an invariant set for a $C^r$-foliation $\F_M$  of a Riemannian manifold $M$, where $r \geq 1$, then the holonomy maps for the   foliation $\F$ on $\fM$   are induced by the holonomy maps of $\F_M$, and there is a metric on the transversals to $\fM$ such that the holonomy maps of $\F$ are Lipschitz, as discussed in  Section~\ref{sec-Lipschitz}.   Problem~\ref{problem7} can be thus be reformulated as.   
 \begin{prob}\label{problem8}
 Let $\fM$ be a Lipschitz matchbox manifold. Find obstructions to the existence of a foliated Lipschitz embedding $\iota \colon \fM \to M$, where $M$ has   a $C^r$-foliation $\F_M$ with $r \geq 1$.     
 \end{prob}
This problem can  also be considered as asking for a  characterization of the Lipschitz structures which can arise for the transverse Cantor sets to exotic minimal sets in $C^r$-foliations.  For example, in the case of a foliation obtained from the suspension of a  diffeomorphism of the circle $\mS^1$,     McDuff studied in \cite{McDuff1981} the question: which Cantor sets embedded in  $\mS^1$ are the invariant sets for  $C^{1+\alpha}$-diffeomorphisms of the circle?  

The general observations and results of this paper are combined  in  Section~\ref{sec-nonembedding} to yield the following non-embedding results.
 
 \begin{thm}\label{thm-noLip1}
There exist  Lipshitz matchbox manifolds which are not homeomorphic to the minimal set of any $C^1$-foliation.
\end{thm}

\begin{thm}\label{thm-noLip2}
There exist  minimal matchbox manifolds which are not homeomorphic to the minimal set of any $C^1$-foliation.
\end{thm}

  Many further questions and problems are posed throughout the text, which is organized as follows.
     
 Section~\ref{sec-foliated}   collects together  some definitions and results concerning matchbox manifolds  that we use in the paper.   More details can be found in   the works \cite{CandelConlon2000,ClarkHurder2013,CHL2013a,CHL2013b,MS2006}. Section~\ref{sec-foliated} is rather dense, and can  be skipped if the reader is only interested in Cantor pseudogroup actions.
 
Section~\ref{sec-dynamics}  gives some some definitions   concerning the dynamical properties of Cantor pseudogroup actions.  Then in  Section~\ref{sec-Lipschitz}    the Lipschitz property for pseudogroup actions is introduced. The main result of this section is a proof that an embedding of a matchbox manifold as an exceptional minimal set in a  $C^1$-foliation yields a Lipschitz structure on it.
 Section~\ref{sec-foliations} discusses some examples from the literature of embeddings of matchbox manifolds as exceptional minimal sets for foliations. 

In Section~\ref{sec-solenoids}, the notion of \emph{normal}, \emph{weak} and \emph{generalized} solenoids are introduced. These are basic examples for the study of minimal matchbox manifolds.

Section~\ref{sec-fusion} introduces an operation on minimal matchbox manifolds, called their ``fusion'', which amalgamates their pseudogroups. The fusion process is inspired by the method introduced by Lukina in \cite{Lukina2012}. 
The fusion process is used to construct the examples in Section~\ref{sec-nonembedding} of minimal pseudogroup   Cantor actions,  which cannot be  homeomorphic to an exceptional minimal set in any $C^1$-foliation.

Finally, in Section~\ref{sec-classification},    \emph{Morita equivalence} and  \emph{Lipschitz equivalence} of minimal Lipschitz pseudogroups are introduced. The problem of the classification of matchbox manifolds up to Lipschitz equivalence is considered   for the the special case of weak solenoids.
 
\section{Foliated spaces and matchbox manifolds} \label{sec-foliated}

We recall the notions of foliated spaces and matchbox manifolds, and their basic properties.
The book   by Moore and Schochet in \cite[Chapter 2]{MS2006} introduced foliated spaces, as part of their development of a general form of the Connes measured leafwise index theorem. The textbook by  Candel and Conlon \cite[Chapter 11]{CandelConlon2000} further develops the theory, with many   examples. Matchbox manifolds are a special class of connected foliated spaces, which have  totally disconnected transversal spaces.
The   papers   \cite{ClarkHurder2011,ClarkHurder2013,CHL2013a,CHL2013b,CHL2013c,CHL2014} discuss the   topology  and dynamics of   \emph{matchbox manifolds},  especially with the goal of classifying these spaces up to homeomorphism.

First we recall some basic notions. A topological  space $\Omega$ is a \emph{continuum},  if it  is \emph{compact, connected, and  metrizable}. 
A Cantor set $\fX$ is a non-empty, compact, perfect and totally disconnected set. A   set $V \subset \fX$ is \emph{clopen} if it is   both open and closed, and a topological space is totally disconnected if and only if it   admits a basis for its topology consisting of clopen sets. 

The definition of a foliated space is modeled on the definition of a smooth foliation.
\begin{defn} \label{def-fs}
A \emph{foliated space of dimension $n$} is a   compact metric space $\fM$, such that  there exists a separable metric space $\fX$, and
for each $x \in \fM$ there is a compact subset $\fX_x \subset \fX$, an open subset $U_x \subset \fM$, and a homeomorphism defined on the closure
$\vp_x \colon \oU_x \to [-1,1]^n \times \fX_x$ such that $\vp_x(x) = (0, w_x)$ where $w_x \in int(\fX_x)$. 
Moreover, it is assumed that each $\vp_x$  admits an extension to a foliated homeomorphism
$\whvarp_x \colon \whU_x \to (-2,2)^n \times \fX_x$ where $\oU_x \subset \whU_x$.
The   space  $\fX_x$ is called   the \emph{local transverse model} at $x$.  
\end{defn}

Let $\pi_x \colon \oU_x \to \fX_x$ denote the composition of $\vp_x$ with projection onto the second factor.
For $w \in \fX_x$ the set $\cP_x(w) = \pi_x^{-1}(w) \subset \oU_x$ is called a \emph{plaque} for the coordinate chart $\vp_x$. We adopt the notation, for $z \in \oU_x$, that $\cP_x(z) = \cP_x(\pi_x(z))$, so that $z \in \cP_x(z)$. Note that each plaque $\cP_x(w)$ for $w \in \fX_x$ is given the topology so that the restriction $\vp_x \colon \cP_x(w) \to [-1,1]^n \times \{w\}$ is a homeomorphism. Then $int (\cP_x(w)) = \vp_x^{-1}((-1,1)^n \times \{w\})$.
Let $U_x = int (\oU_x) = \vp_x^{-1}((-1,1)^n \times int(\fX_x))$.
Note that if $z \in U_x \cap U_y$, then $int(\cP_x(z)) \cap int( \cP_y(z))$ is an open subset of both
$\cP_x(z) $ and $\cP_y(z)$.
The collection of sets
$$\cV = \{ \vp_x^{-1}(V \times \{w\}) \mid x \in \fM ~, ~ w \in \fX_x ~, ~ V \subset (-1,1)^n ~ {\rm open}\}$$
forms the basis for the \emph{fine topology} of $\fM$. The connected components of the fine topology are called \emph{leaves}, and define the foliation $\F$ of $\fM$.
Let $L_x \subset \fM$ denote the leaf of $\F$ containing $x \in \fM$.

\begin{defn} \label{def-sfs}
A \emph{smooth foliated space} is a foliated space $\fM$ as above, for which there exists a choice of local charts $\vp_x \colon \oU_x \to [-1,1]^n \times \fX_x$ such that for all $x,y \in \fM$ with $z \in U_x \cap U_y$, there exists an open set $z \in V_z \subset U_x \cap U_y$ such that $\cP_x(z) \cap V_z$ and $\cP_y(z) \cap V_z$ are connected open sets, and the composition
$\ds \psi_{x,y;z} \equiv \vp_y \circ \vp_x ^{-1}\colon \vp_x(\cP_x (z) \cap V_z) \to \vp_y(\cP_y (z) \cap V_z)$
is a smooth map, where $\vp_x(\cP_x (z) \cap V_z) \subset \mR^n \times \{w\} \cong \mR^n$ and $\vp_y(\cP_y (z) \cap V_z) \subset \mR^n \times \{w'\} \cong \mR^n$. The  maps $\psi_{x,y;z}$ are assumed to depend continuously on $z$ in the $C^{\infty}$-topology on maps between subsets of $\mR^n$.
\end{defn}

A map $f \colon \fM \to \mR$ is said to be \emph{smooth} if for each flow box
$\vp_x \colon \oU_x \to [-1,1]^n \times \fX_x$ and $w \in \fX_x$ the composition
$y \mapsto f \circ \vp_x^{-1}(y, w)$ is a smooth function of $y \in (-1,1)^n$, and depends continuously on $w$ in the $C^{\infty}$-topology on maps of the plaque coordinates $y$. As noted in \cite{MS2006} and \cite[Chapter 11]{CandelConlon2000}, this allows one to define smooth partitions of unity, vector bundles, and tensors for smooth foliated spaces. In particular, one can define leafwise Riemannian metrics. We recall a standard result, whose proof for foliated spaces can be found in \cite[Theorem~11.4.3]{CandelConlon2000}.
\begin{thm}\label{thm-riemannian}
Let $\fM$ be a smooth foliated space. Then there exists a leafwise Riemannian metric for $\F$, such that for each $x \in \fM$, $L_x$ inherits the structure of a complete Riemannian manifold with bounded geometry, and the Riemannian geometry of $L_x$ depends continuously on $x$. In particular, each leaf $L_x$ has the structure of a complete Riemannian manifold with bounded geometry. 
\end{thm}

Bounded geometry implies, for example, that for each $x \in \fM$, there is a leafwise exponential map
$\exp^{\F}_x \colon T_x\F \to L_x$ which is a surjection, and the composition $\exp^{\F}_x \colon T_x\F \to L_x \subset \fM$ depends continuously on $x$ in the compact-open topology on maps.

\begin{defn} \label{def-mm}
A \emph{matchbox manifold} is a  
smooth foliated connected space $\fM$, such that  its  transverse model space $\fX$ is totally disconnected, and for each $x \in \fM$, the transverse model space $\fX_x \subset \fX$ in Definition~\ref{def-fs} is a clopen  subset.
\end{defn}

All matchbox manifolds  are assumed to be smooth with a given leafwise Riemannian metric. 
The space $\fM$ is assumed to be metrizable, and we fix a  choice  for the metric  $\dM$ on $\fM$. 
The leafwise Riemannian metric $\dF$ is continuous with respect to  the metric $\dM$ on $\fM$, but otherwise    the two metrics  can be chosen independently. 
 The metric $\dM$ is       used to define the metric topology on $\fM$, while the metric $\dF$ depends on an independent  choice of the Riemannian metric on leaves.

An important difference between a foliated matchbox manifold and a smooth foliated manifold is that the local foliation charts for a matchbox manifold are not connected, and so must be chosen appropriately to ensure that each chart is ``local''.   We introduce the following conventions. 
 
For $x \in \fM$ and $\e > 0$, let $D_{\fM}(x, \e) = \{ y \in \fM \mid \dM(x, y) \leq \e\}$ be the closed $\e$-ball about $x$ in $\fM$, and $B_{\fM}(x, \e) = \{ y \in \fM \mid \dM(x, y) < \e\}$ the open $\e$-ball about $x$.

Similarly, for $w \in \fX$ and $\e > 0$, let $D_{\fX}(w, \e) = \{ w' \in \fX \mid d_{\fX}(w, w') \leq \e\}$ be the closed $\e$-ball about $w$ in $\fX$, and $B_{\fX}(w, \e) = \{ w' \in \fX \mid d_{\fX}(w, w') < \e\}$ the open $\e$-ball about $w$.

Given a leaf $L$ and a piecewise $C^1$-path $\gamma \colon [0,1] \to L$, let  $\| \gamma \|_{\F}$ denote its  path-length for the leafwise Riemannian metric.  Then  give $L \subset \fM$ the     path-length metric: if $x, y \in L$ then set
$$\dF(x,y) = \inf \left\{\| \gamma\|_{\F} \mid \gamma \colon [0,1] \to L ~{\rm is ~ piecewise ~~ C^1}~, ~ \gamma(0) = x ~, ~ \gamma(1) = y ~, ~ \gamma(t) \in L \quad \forall ~ 0 \leq t \leq 1\right\},$$
and otherwise, if  $x,y \in \fM$   are not on the same leaf, then set $\dF(x,y) = \infty$. 
 
   For each $x \in \fM$ and $r > 0$, let $D_{\F}(x, r) = \{y \in L_x \mid \dF(x,y) \leq r\}$.

For each $x \in \fM$, the  {Gauss Lemma} implies that there exists $\lambda_x > 0$ such that $D_{\F}(x, \lambda_x)$ is a \emph{strongly convex} subset for the metric $\dF$. That is, for any pair of points $y,y' \in D_{\dF}(x, \lambda_x)$ there is a unique shortest geodesic segment in $L_x$ joining $y$ and $y'$ and  contained in $D_{\F}(x, \lambda_x)$ (cf. \cite[Chapter 3, Proposition 4.2]{doCarmo1992}, or \cite[Theorem 9.9]{Helgason1978}). Then for all $0 < \lambda < \lambda_x$ the disk $D_{\F}(x, \lambda)$ is also strongly convex. The leafwise metrics have uniformly bounded geometry, so we obtain:

\begin{lemma}\label{lem-stronglyconvex}
There exists $\lF > 0$ such that for all $x \in \fM$, $D_{\F}(x, \lF)$ is strongly convex.  
\end{lemma}

The following proposition summarizes results in \cite[sections 2.1 - 2.2]{ClarkHurder2013}.

\begin{prop}\label{prop-regular} 
For a smooth foliated space $\fM$, given $\eM > 0$, there exist  constants  $\lF>0$ and $0< \dFU < \lF/5$, and  a covering of $\fM$ by foliation  charts 
$\ds \left\{\vp_i \colon \oU_i \to [-1,1]^n \times \fX_i \mid 1 \leq i \leq \nu \right\}$
 with the following properties: For each  $1 \leq i \leq \nu$, let   $\pi_i = \pi_{x_i} \colon \oU_i \to \fX_i$ be the projection, then
\begin{enumerate}
\item  Interior: $U_i \equiv int(\oU_i) = \vp_i^{-1}\left( (-1,1)^n \times B_{\fX}(w_i, \e_i)\right)$, where  $w_i \in \fX_i$ and $\e_i>0$.
\item Locality: for   $x_i \equiv \vp_i^{-1}(w_i, 0) \in \fM$,    $\oU_i \subset B_{\fM}(x_i, \eM)$.
\end{enumerate}
For $z \in \oU_i$, the \emph{plaque} of the chart $\vp_i$ through $z$ is denoted by $\cP_i(z) = \cP_i(\pi_i(z)) \subset \oU_i$. 
\begin{enumerate}\setcounter{enumi}{2}
\item Convexity:  the plaques of $\vp_i$ are   strongly convex subsets for the leafwise metric. 
\item Uniformity:  for   $w \in \fX_i$ let $x_{w} = \vp_{x_i}^{-1}(0 , w)$, then 
\begin{equation}\label{eq-Fdelta}
D_{\F}(x_{w} , \dFU/2) ~ \subset ~ \cP_i(w) ~ \subset ~ D_{\F}(x_{w} , \dFU) 
\end{equation} 
\item \label{item-clopen} The projection $\pi_i(U_i \cap U_j) = \fX_{i,j} \subset \fX_i$ is a clopen subset for all $1 \leq i, j \leq \nu$. 
\end{enumerate}
A \emph{regular foliated covering} of $\fM$ is one that satisfies   the above conditions (\ref{prop-regular}.1) to (\ref{prop-regular}.5). 
\end{prop}
This technical result highlights one of the main issues with foliated spaces and matchbox manifolds, as contrasted with smooth foliations of compact manifolds,   one has to assume or prove for these more exotic foliated spaces many of the regularity properties that are used in the study of foliations.

We assume in the following that a   regular foliated covering  of $\fM$ as in Proposition~\ref{prop-regular}  has been chosen.
 Let $\cU = \{U_{1}, \ldots , U_{\nu}\}$ denote the corresponding open covering of $\fM$.
We can   assume that the spaces $\fX_i$ form a \emph{disjoint  clopen covering} of $\fX$, so that 
  $\ds \fX = \fX_1 \dot{\cup} \cdots \dot{\cup} \fX_{\nu}$.

Let $\eU > 0$ be a Lebesgue number for $\cU$. That is, given any $z \in \fM$ there exists some index $1 \leq i_z \leq \nu$ such that the open metric ball $B_{\fM}(z, \eU) \subset U_{i_z}$.

For $1 \leq i \leq \nu$, let  $ \lambda_i \colon \oU_i \to [-1,1]^n$ be the projection, so that for each $z \in U_i$   the restriction $\lambda_i \colon \cP_i(z) \to [-1,1]^n$ is  is a smooth coordinate system on the plaque.

For each $1 \leq i \leq \nu$ the set $\cT_i =  \vp_i^{-1}(0 , \fX_i)$ is a compact transversal to $\F$. Without loss of generality, we can assume   that the transversals 
$\ds \{ \cT_{1} , \ldots , \cT_{\nu} \}$ are pairwise disjoint in $\fM$. Then define sections
\begin{equation}\label{eq-taui}
\tau_i \colon \fX_i \to \oU_i ~ , ~ {\rm defined ~ by} ~ \tau_i(\xi) = \vp_i^{-1}(0 , \xi) ~ , ~ {\rm so ~ that} ~ \pi_i(\tau_i(\xi)) = \xi.
\end{equation}
Then $\cT_i = \cT_{x_i}$ is the image of $\tau_i$ and we let $\cT = \cT_1 \cup \cdots \cup \cT_{\nu} \subset \fM$ denote their disjoint union, and $\tau \colon \fX \to \cT$ the union of the maps $\tau_i$.
   
A map $f \colon \fM \to \fM'$ between foliated spaces is said to be a \emph{foliated map} if the image of each leaf of $\F$ is contained in a leaf of $\F'$. If $\fM'$ is a matchbox manifold, then each leaf of $\F$ is path connected, so its image is path connected, hence must be contained in a leaf of $\F'$. Thus, 
\begin{lemma} \label{lem-foliated1}
Let $\fM$ and $\fM'$ be matchbox manifolds, and $f \colon \fM' \to \fM$ a continuous map. Then $f$ maps the leaves of $\F'$ to leaves of $\F$. In  particular, any homeomorphism $f \colon \fM' \to \fM$ of   matchbox manifolds is a foliated map. \hfill $\Box$
\end{lemma}

A \emph{leafwise path}  is a continuous map $\gamma \colon [0,1] \to \fM$ such that there is a leaf $L$ of $\F$ for which $\gamma(t) \in L$ for all $0 \leq t \leq 1$. 
If $\fM$ is a matchbox manifold, and $\gamma \colon [0,1] \to \fM$ is continuous, then   $\gamma$ is a leafwise path by Lemma~\ref{lem-foliated1}. In the following, we will assume that all paths are  piecewise differentiable.

The holonomy pseudogroup of a smooth foliated manifold $(M, \F)$ generalizes the concept of a Poincar\'{e}   section for a flow, which induces a discrete  dynamical system    associated to the flow. Associated  to a leafwise path $\gamma$ is a holonomy map $h_{\gamma}$, which is a local homeomorphism  on the transversal space. For a matchbox manifold $(\fM, \F)$ the holonomy along a leafwise path  is defined analogously.  We briefly recall below the ideas and notations of the construction  of holonomy maps for matchbox manifolds;   further details  and  proofs are given in  \cite{ClarkHurder2013,CHL2013a}.

 A pair of indices $(i,j)$, $1 \leq i,j \leq \nu$, is said to be \emph{admissible} if  $U_i \cap U_j \ne \emptyset$.
For $(i,j)$ admissible, set $\fX_{i,j} = \pi_i(U_i \cap U_j) \subset \fX_i$.  The regularity of foliation charts imply that plaques are either disjoint, or have connected intersection. For $(i,j)$ admissible,   there is a well-defined transverse change of coordinates homeomorphism $h_{i,j} \colon \fX_{i,j} \to \fX_{j,i}$ with domain $\Dom(h_{i,j}) = \fX_{i,j}$ and range $R(h_{i,j}) = \Dom(h_{j,i}) = \fX_{j,i}$.   By definition they satisfy $h_{i,i} = Id$, $h_{i,j}^{-1} = h_{j,i}$, and if $U_i \cap U_j\cap U_k \ne \emptyset$ then $h_{k,j} \circ h_{j,i} = h_{k,i}$ on their common domain of definition. 
Note that the domain and range of $h_{i,j}$ are clopen subsets of $\fX$ by Proposition~\ref{prop-regular}.\ref{item-clopen}.

 Recall that for $1 \leq i \leq \nu$, $\tau_i \colon \fX_i \to \cT_i$ denotes the transverse section for the coordinate chart $U_i$, where  $\cT = \cT_1 \cup \cdots \cup \cT_{\nu} \subset \fM$ denotes their disjoint union, and  $\pi \colon \cT \to \fX$ is the coordinate projection restricted to $\cT$ which is a homeomorphism, with $\tau \colon \fX \to \cT$   its inverse.

The \emph{holonomy pseudogroup} $\cGF$ of $\F$ is the topological pseudogroup   modeled on $\fX$ generated by   the elements of   $\cGF^{(1)} = \{h_{j,i} \mid (i,j) ~{\rm admissible}\}$.  We also define a  subpseudogroup $\ds \cGF^*  \subset \cGF$ which is based on the holonomy along paths. 
A sequence $\cI = (i_0, i_1, \ldots , i_{\alpha})$ is \emph{admissible} if each pair $(i_{\ell -1}, i_{\ell})$ is admissible  for $1 \leq \ell \leq \alpha$, and the composition
$\ds  h_{\cI} = h_{i_{\alpha}, i_{\alpha-1}} \circ \cdots \circ h_{i_1, i_0}$ 
 has non-empty domain $\Dom(h_{\cI})$, which is defined to be    the  maximal clopen subset of $\fX_{i_0}$ for which the compositions are defined.
Given a  open subset $U \subset \Dom(h_{\cI})$ define the restriction $h_{\cI} | U \in \cGF$.  Introduce
\begin{equation}\label{eq-restrictedpseudogroup}
\cGF^* = \left\{ h_{\cI} |  U \mid   \cI ~ {\rm admissible~ and} ~ U \subset \Dom(h_{\cI}) \right\} \subset \cGF ~ .
\end{equation}
 The range of $g = h_{\cI} |  U$ is the open set $R(g) = h_{\cI}(U) \subset \fX_{i_{\alpha}} \subset \fX$. Note that each map $g \in \cGF^*$ admits a
continuous extension $\overline{g} \colon \overline{\Dom(g)} = \overline{U} \to \fX_{i_{\alpha}}$ as $\Dom( h_{\cI})$ is a clopen set for each $\cI$.

Let $\cI = (i_0, i_1, \ldots , i_{\alpha})$ be an  admissible sequence. 
For each $1 \leq \ell \leq \alpha$, set 
$\cI_{\ell} = (i_0, i_1, \ldots, i_{\ell})$, and let $h_{\cI_{\ell}}$ denote the corresponding holonomy map. For $\ell = 0$, let $\cI_0 = (i_0 , i_0)$.
Note that $h_{\cI_{\alpha}} = h_{\cI}$ and $h_{\cI_{0}} = Id \colon \fX_0 \to \fX_0$.

Given $w \in \Dom(h_{\cI})$,  let $x = \tau_{i_0}(w) \in L_{w}$. For each 
$0 \leq \ell \leq \alpha$, set $w_{\ell} = h_{\cI_{\ell}}(w)$ and
$x_{\ell}= \tau_{i_{\ell}}(w_{\ell})$. 
Recall that $\cP_{i_{\ell}}(x_{\ell}) = \cP_{i_{\ell}}(w_{\ell})$, where 
each $\cP_{i_{\ell}}(w_{\ell})$ is a strongly convex subset of the   leaf $L_w$ in the leafwise metric $d_{\F}$. 
 Introduce the   \emph{plaque chain}
 \begin{equation}\label{eq-plaquechain}
\cP_{\cI}(w) = \{\cP_{i_0}(w_0), \cP_{i_1}(w_1), \ldots , \cP_{i_{\alpha}}(w_{\alpha}) \} ~ .
\end{equation}
Adopt the notation $\cP_{\cI}(x) \equiv \cP_{\cI}(w)$.
  Intuitively, a plaque chain $\cP_{\cI}(x)$ is a sequence of successively overlapping convex ``tiles'' in $L_{w}$ starting at $x = \tau_{i_0}(w)$, ending at
$y = x_{\alpha} = \tau_{i_{\alpha}}(w_{\alpha})$, and with each $\cP_{i_{\ell}}(x_{\ell})$ ``centered'' on the point $x_{\ell} = \tau_{i_{\ell}}(w_{\ell})$.

    Let $\gamma \colon  [0,1] \to \fM$ be a   path. Set $x_0 = \gamma(0) \in U_{i_0}$,   $w = \pi(x_0)$ and $x = \tau(w) \in \cT_{i_0}$.  
Let $\cI$ be an admissible sequence with  $w \in \Dom(h_{\cI})$. We say that $(\cI , w)$ \emph{covers} $\gamma$ 
if the domain of $\gamma$ admits   a partition $0 = s_0 < s_1 < \cdots < s_{\alpha} = 1$ such that   $\cP_{\cI}(w)$   satisfies
\begin{equation}\label{eq-cover}
\gamma([s_{\ell} , s_{\ell + 1}]) \subset    \cP_{i_{\ell}}(\xi_{\ell})  ~ , ~ 0 \leq \ell < \alpha, ~ {\rm and} ~  \gamma(1) \in   \cP_{i_{\alpha}}(\xi_{\alpha}) .
\end{equation}

  For  a   path $\gamma$, we   construct an admissible sequence
$\cI = (i_0, i_1, \ldots, i_{\alpha})$ with $w \in \Dom(h_{\cI})$ so that $(\cI , w)$ covers $\gamma$, and has ``uniform domains''.
Inductively choose a partition of the interval $[0,1]$, say $0 = s_0 < s_1 < \cdots < s_{\alpha} = 1$,  such that for each $0 \leq \ell \leq \alpha$,
$$\gamma([s_{\ell}, s_{\ell + 1}]) \subset D_{\F}(x_{\ell}, \eFU) \quad , \quad x_{\ell} = \gamma(s_{\ell}).$$
As a notational convenience, we have let
$s_{\alpha+1} = s_{\alpha}$, so that $\gamma([s_{\alpha}, s_{\alpha + 1}]) = x_{\alpha}$.
Choose $s_{\ell + 1}$ to be the largest value of $s_{\ell} < s \leq 1$ such that $\dF(\gamma(s_{\ell}), \gamma(t)) \leq \eFU$ for all  $s_{\ell} \leq t \leq s$,  then    $\alpha \leq   \| \gamma \|/\eFU$. 

For each $0 \leq \ell \leq \alpha$, choose an index $1 \leq i_{\ell} \leq \nu$ so that $ B_{\fM}(x_{\ell}, \eU) \subset U_{i_{\ell}}$.
Note that, for all $s_{\ell} \leq t \leq s_{\ell +1}$, $B_{\fM}(\gamma(t), \eU/2) \subset U_{i_{\ell}}$, so that
$x_{\ell+1} \in U_{i_{\ell}} \cap U_{i_{\ell +1}}$. It follows that $\cI_{\gamma} = (i_0, i_1, \ldots, i_{\alpha})$ is an admissible sequence.
Set $h_{\gamma} = h_{\cI_{\gamma}}$ and note that  $h_{\gamma}(w) = w'$.

Next, consider paths $\gamma, \gamma' \colon [0,1] \to \fM$   with $x = \gamma(0) = \gamma'(0)$ and $y = \gamma(1) = \gamma'(1)$. Suppose that $\gamma$ and $\gamma'$ are homotopic relative endpoints. That is, assume  there exists a continuous map $H \colon [0,1] \times [0,1] \to \fM$ with 
$$H(0,t) = \gamma(t) ~, ~ H(1,t) = \gamma'(t) ~ , ~ H(s,0) = x ~ {\rm and} ~ H(s,1) = y \quad {\rm for ~ all} ~ 0 \leq s \leq 1$$
Then there exists partitions $0 = s_0 < s_1 < \cdots < s_{\beta} = 1$ and $0 = t_0 < t_1 < \cdots < t_{\alpha} = 1$ such that for each pair of indices $0 \leq j < \beta$ and $0 \leq k < \alpha$,   there is an index $1 \leq i(j,k)\leq \nu$ such that 
$$H([s_j,s_{j+1}] \times [t_k, t_{k+1}] ) \subset D_{\F}(H(s_j, t_k), \eFU) \subset U_{i(j,k)}$$
A standard argument then yields     the following basic fact about holonomy maps. 
\begin{lemma}\label{lem-homotopy}
Let $\gamma, \gamma' \colon [0,1] \to \fM$ be   paths with $x = \gamma(0) = \gamma'(0)$ and $y = \gamma(1) = \gamma'(1)$, and suppose they are homotopic relative endpoints. Then the induced holonomy maps $h_{\gamma}$ and $h_{\gamma'}$ agree on an open neighborhood of $\xi_0 = \pi_{i_0}(x)$. 
\end{lemma}

Next consider the \emph{groupoid} formed by germs of maps in $\cGF$.  
  Let  $U, U', V, V' \subset \fX$ be open subsets with  $w \in U \cap U'$. Given homeomorphisms    $h \colon U \to V$  and $h' \colon U' \to V'$    with $h(w) = h'(w)$,   then   $h$ and $h'$ have the same \emph{germ at $w$}, and write    $h \sim_w h'$,   if there exists an open neighborhood $w \in W \subset U \cap U'$ such that $h | W= h' |W$. Note that $\sim_w$ defines an equivalence relation. 

\begin{defn}\label{def-germ}
The \emph{germ of $h$ at $w$} is the equivalence class $[h]_w$ under the relation ~$\sim_w$. The  map  $h \colon U \to V$ is called a \emph{representative} of  $[h]_w$.
The point $w$ is called the source of  $[h]_w$ and denoted $s([h]_w)$, while $w' = h(w)$ is called the range of  $[h]_w$ and denoted $r([h]_w)$.
\end{defn}

   The collection of all such germs $[h]_w$ for   $h \in \cGF$ and $w \in \Dom(h)$,  forms  the \emph{holonomy   groupoid} $\GF$, which has the natural topology associated to sheaves of maps over $\cX$. Let $\cRF \subset \fX \times \fX$ denote the equivalence relation on  $\fX$ induced by $\F$, where     $(w,w') \in \cRF$ if and only if $w,w'$ correspond to points on the same leaf of $\F$. The product map $s \times r \colon \GF \to \cRF$ is  \'etale; that is, the map is a local homeomorphism with discrete fibers.
These notions were  introduced by Haefliger for foliations  \cite{Haefliger1958,Haefliger1984}, and naturally extend to the case of matchbox manifolds.
 
 We introduce a convenient notation for elements of $\GF$.
 Let $(w,w') \in \cRF$, and let $\gamma$  denote a path from $x = \tau(w)$ to $y = \tau(w')$. We may assume that $\gamma$ is a geodesic for the leafwise metric, 
 and let $[h_{\gamma}]_w$ (or sometimes just $\gamma_w$) denote the germ at $w$ of the holonomy map defined by $\gamma$. 
 
It follows that  there is a well-defined surjective homomorphism, the \emph{holonomy map},  
 \begin{equation}\label{eq-holodef}
h_{\F,x} \colon \pi_1(L_x , x) \to \G_w^w \equiv \left\{  [g]_w \in \GF  \mid    r([g]_w) =w \right\}   
\end{equation}
Moreover,  if $y,z \in L$ then the homomorphism
$h_{\F , y}$ is conjugate (by an element of $\cGF$) to the homomorphism $h_{\F , z}$.
A leaf $L$ is said to have \emph{non-trivial germinal holonomy} if for some $y \in L$, the homomorphism $h_{\F , y}$ is non-trivial. If the homomorphism $h_{\F , y}$ is trivial, then we say that $L_y$ is a \emph{leaf without holonomy}. This property depends only on $L$, and not the choice of  $y \in L$.

 \begin{lemma}\label{lem-homotopymin}
Given a path $\gamma \colon [0,1] \to \fM$  with $x = \gamma(0)$ and $y = \gamma(1)$. Suppose that   $L_x$ is a leaf without holonomy. Then there exists a leafwise geodesic segment $\gamma'  \colon [0,1] \to \fM$  with $x = \gamma'(0)$ and $y = \gamma'(1)$,  such that $\|\gamma' \| = \dF(x,y)$,    and $h_{\gamma}$ and $h_{\gamma'}$ agree on an open neighborhood of $\xi_0$.
\end{lemma}
\proof
The leaf $L_x$ containing $x$ is a complete Riemannian manifold, so there exists  a geodesic segment $\gamma'$ which is length minimizing between $x$ and $y$.
Then the holonomy maps  $h_{\gamma}$ and $h_{\gamma'}$ agree on an open neighborhood of $\xi_0  = \pi_{i_0}(x)$ by the definition of germinal holonomy. 
\endproof

   Next, we introduce the filtrations of $\cGF^*$   by word length, and of $\GF$ by    path length, then derive estimates comparing these notions of length.

For $\alpha \geq 1$, let  $\cGF^{(\alpha)}$ be the collection of holonomy homeomorphisms $h_{\cI} | U \in \cGF^*$ determined   by admissible paths $\cI = (i_0,\ldots,i_k)$ such that $k \leq \alpha$ and $U \subset \Dom(h_{\cI})$ is open. For each $\alpha$, let $C(\alpha)$ denote the number of admissible sequences of length at most $\alpha$. As there are at most $\nu^2$ admissible pairs $(i,j)$, we have the basic estimate that $C(\alpha) \leq \nu^{2 \alpha}$. This upper bound estimate grows exponentially with $\alpha$, though the exact growth rate of $C(\alpha)$ may be much less.

For each $g \in \cGF^*$ there is some $\alpha$ such that $g \in \cGF^{(\alpha)}$. Let $\|g\|$ denote the least such $\alpha$, which is called the \emph{word length} of $g$.  Note that  $\cGF^{(1)}$ generates $\cGF^*$.

We use the word length on $\cGF^*$ to  define the word length on  $\GF$, where for  $\gamma_w \in \GF$, set
\begin{equation}
\| \gamma_w \| ~ = ~ \min ~ \left\{ \| g \| \mid [g]_w = \gamma_w ~ {\rm for}~ g \in \cGF^*  \right\} .
\end{equation}

  Introduce the \emph{path length} of $\gamma_w \in \GF$, by considering the infimum of the lengths $\| \gamma'\|$ for all   piecewise smooth curves $\gamma'$  for which $\gamma_w' = \gamma_w$. That is, 
\begin{equation}\label{eq-groupoidpathlength}
\ell(\gamma_w)  ~ = ~ \inf ~ \left\{ \| \gamma' \| \mid \gamma'_w = \gamma_w   \right\} .
\end{equation}
Note that if $L_w$ is a leaf without holonomy, set $x = \tau(w)$ and $y = \tau(w')$, then Lemma~\ref{lem-homotopymin} implies that 
$\ell(\gamma_w) = \dF(x,y)$. 
 This yields a   fundamental estimate, whose proof can be found in \cite{CHL2013b}:
\begin{lemma}\label{lem-comparisons}
Let  $[g]_w \in \GF$ where $w$ corresponds to a leaf without holonomy. Then
\begin{equation}\label{eq-comparisons}
\dF(x,y)/2\dFU ~ \leq ~ \| [g]_w \| ~ \leq ~  1 + \dF(x,y)/\eFU
\end{equation}
\end{lemma}

\section{Pseudogroup Dynamics} \label{sec-dynamics}

In this section, we consider some   aspects of the topological dynamics of pseudogroups, which are useful for obtaining dynamical invariants for   the pseudogroup $\cGX$  associated to a matchbox manifold. 
The sources \cite{CandelConlon2000,Hurder2014,Walczak2004} give more detailed discussions. 
 
  The study of the dynamics of a pseudogroup $\cGX$ acting on $\fX$  is  a generalization of the study of continuous actions of finitely-generated groups on Cantor sets, though it differs   in some fundamental ways.  For a group action, each $\gamma \in \G$ defines a homeomorphism $h_{\gamma} \colon \fX \to \fX$.  For a pseudogroup action,  given $g  \in \cGX$ and $w \in \Dom(g)$, 
 there is some clopen neighborhood $w \in U \subset \Dom(g)$ for which $g | U = h_{\cI} | U$ where $\cI$ is   admissible sequence  with $w \in \Dom(h_{\cI})$. 
 By the definition of a pseudogroup, every $g \in \cGX$ is the ``union'' of such maps, and the dynamical properties of the action may reflect the fact that the domains of the actions are not all of $\fX$.

  We first recall some basic definitions.
 
 \begin{defn}\label{def-lipequiv}
 A pseudogroup $\cGX$ acting on a Cantor set $\fX$ is  \emph{compactly generated}, if there exists two collections of \emph{clopen} subsets,  $\{U_1, \ldots, U_k\}$ and $\{V_1, \ldots, V_k\}$ of $\fX$, and homeomorphisms $\{h_i \colon U_i \to V_i \mid 1 \leq i \leq k\}$ which generate all elements of $\cGX$.
The collection of maps $\cGX^*$  is defined to be all  compositions of the generators on the maximal domains for which the composition is defined. 
 \end{defn}

Let $\dX$ be a metric on $\fX$ which defines the topology on the space.

\begin{defn} \label{def-expansive}
The action of a compactly generated pseudogroup  $\cGX$ on $\fX$ is \emph{expansive}, or more properly \emph{$\e$-expansive}, if there exists $\e > 0$ such that for all $w, w' \in \fX$, there exists $g \in \cGX^*$ with $w, w' \in D(g)$ such that $\dX(g(w), g(w')) \geq \e$.
\end{defn}

\smallskip

\begin{defn} \label{def-equicontinuous}
The action of a compactly generated pseudogroup  $\cGX$ on $\fX$  is  \emph{equicontinuous} if for all $\e > 0$, there exists $\delta > 0$ such that for all $g \in \cGX^*$, if $w, w' \in D(g)$ and $\dX(w,w') < \delta$, then $\dX(g(w), g(w')) < \e$.
Thus, $\cGX^*$ is equicontinuous as a family of local group actions. 
\end{defn}

The  \emph{geometric entropy}   for  pseudogroup actions, introduced by Ghys, Langevin and Walczak \cite{GLW1988},  gives a measure    of the ``exponential complexity'' of the orbits of the action. See also the discussion of entropy for pseudogroup actions  in  Candel and Conlon \cite[\S 13.2B]{CandelConlon2000}, and in  Walczak \cite{Walczak2004}.

The key idea is the notion of $\e$-separated sets, due to Bowen \cite{Bowen1971}. 
  Let $\e > 0$ and $\ell > 0$. 
   A subset $\cE \subset \fX$ is said to be $(\dX, \e, \ell)$-separated if   for all $w,w' \in \cE \cap \fX_i$   there exists $g \in \cGX^*$ with $w,w' \in \mathrm{Dom}(g) \subset \fX_i$, and  $\|g\|_w \leq \ell $  so that  $\dX(g(w), g(w')) \geq \e$. 
  If $w  \in \fX_i$ and $w' \in \fX_j$ for $i \ne j$ then   they are   $(\e, \ell)$-separated by default.   
The  ``expansion growth function'' counts the maximum  of this quantity:
  $$h(\cGX, \dX, \e, \ell) =     \max \{ \# \cE \mid \cE \subset \fX ~ \text{is} ~ (\dX, \e,\ell) \text{-separated} \}   $$
The entropy is then defined to be  the exponential growth type of the expansion growth function:
$$ h(\cGX, \dX, \e)   =   \limsup_{\ell \to  \infty} ~  \ln \left\{ h(\cGX, \dX, \e, \ell) \right\}/ \ell    \quad , \quad h(\cGX,\dX)   =  \lim_{\e \to 0} ~  h(\cGX, \dX, \e) $$
Note that the quantity $h(\cGX, \dX) \geq 0$, and it  may take the value $h(\cGX, \dX)  = \infty$. 

We recall two key properties of pseudogroup entropy. The first property follows directly from the definition of entropy.

\begin{prop}[Proposition~2.6, \cite{GLW1988}]\label{prop-metricentropy1}
Let    $\cGX$ be a compactly generated pseudogroup, acting on the compact space  $\fX$ with the metric $\dX$. Then the geometric entropy $h(\cGX,\dX)$ is independent of the choice of metric $\dX$. 
\end{prop}

The second property is an exercise using standard properties of the  pseudogroup length function.
\begin{prop}[Exercise~13.2.21, \cite{CandelConlon2000}]\label{prop-metricentropy2}
Let    $\cGX$ be a compactly generated pseudogroup, acting on a compact space  $\fX$ with the metric $\dX$. Then the property that $h(\cGX,\dX)$ is either zero, finite, or infinite, is independent of the choice of   generating set for   $\cGX$. 
\end{prop}

\section{Lipschitz foliations and geometry} \label{sec-Lipschitz}

In this section, we define the Lipschitz property for   matchbox manifolds $\fM$. The basic result is that if $\fM$ is homeomorphic to an exceptional minimal set in a $C^1$-foliation, then its transversal space $\fX$ has a metric for which the induced pseudogroup  $\cG_\fX$ is Lipschitz.  

It is a standard fact that there is a  {unique} Cantor set, up to \emph{homeomorphism}. That is, any two compact, perfect, totally disconnected and non-empty sets are homeomorphic. See \cite[Chapter~12]{Moise1977} for a proof and discussion of this result.  In particular, for a given Cantor set $\fX$, any non-empty clopen  subset $U \subset \fX$ is homeomorphic to $\fX$. 

Two metrics $\dX$ and $\dXp$ are  \emph{Lipschitz equivalent}  if     for some   $C \geq 1$,
they satisfy the condition:
\begin{equation}
C^{-1} \cdot \dX(x,y) ~ \leq ~ \dXp(x,y) ~ \leq ~ C \cdot \dX(x,y) \quad {\rm for ~ all} ~ x,y \in \fX
\end{equation}

On the other hand, there are many possible metrics on $\fX$ which are compatible with its topology, and need not be Lipschitz equivalent. 
The study of the \emph{Lipschitz geometry} of the pair $(\fX, \dX)$  investigates the geometric properties common to all metrics in the Lipschitz class of the given metric $\dX$.  Problem~\ref{problem8} can be rephrased as asking for characterizations of the transverse Lipschitz geometry of exceptional minimal sets.

 We next consider the Lipschitz property of matchbox manifolds. 
  The choice of a regular foliated covering  
$\ds \left\{\vp_i \colon \oU_i \to [-1,1]^n \times \fX_i \mid 1 \leq i \leq \nu \right\}$
for the   matchbox manifold $\fM$, as in Proposition~\ref{prop-regular}, yields the pseudogroup $\cGX$ which acts via homeomorphisms on the transversal space $\fX$ to $\F$.

\begin{defn} \label{def-Lipschitz}
 The action of a  compactly generated pseudogroup  $\cGX$  is \emph{C-Lipschitz}  with respect to $\dX$, if there exists a generating set 
 $\{h_i \colon U_i \to V_i \mid 1 \leq i \leq k\}$  as in Definition~\ref{def-lipequiv}, and $C \geq 1$, such that for each   $1 \leq i \leq k$  and for all $w, w' \in U_i$ we have
\begin{equation}\label{eq-Lipschitz}
C^{-1} \cdot \dX(w,w') \leq \dX(h_i(w), h_i(w')) \leq C \cdot \dX(w,w') ~ .
\end{equation}
\end{defn}
The condition \eqref{eq-Lipschitz} is equivalent to saying that $\cGX$ is generated by \emph{bi-Lipschitz homeomorphisms}, though we use the notation  Lipschitz for the action of the pseudogroup $\cGX$.

 Recall that $\tau \colon \fX \to \cT \subset \fM$ is the transversal to $\F$ associated to a regular covering of $\fM$.
 Let $\dX$ be the metric induced on $\fX$ by the restriction of $\dM$ on $\fM$ to the image of     $\tau \colon \fX \to \cT$.

  The claim of the following result is intuitively clear, but its proof requires   care with the subtleties of working with foliation charts that have totally disconnected transversals.

 \begin{prop}\label{prop-lipembed}
Let $\fM$ be a minimal matchbox manifold, and $M$ a smooth Riemannian manifold with a $C^1$-foliation $\F_M$, and $\cZ \subset M$ an exceptional minimal set for $\F_M$.
Suppose there exists  a homeomorphism $f \colon \fM \to \cZ$, then there exists a metric $\dX$ on $\fX$ 
such that the action of  the holonomy pseudogroup  $\cGX$ on $\fX$    is Lipschitz.
 \end{prop}
\proof
The map $f$   maps leaves of $\F$ to leaves of $\F_M$ by Lemma~\ref{lem-foliated1}, and as $f$ is a homeomorphism onto its image, this  implies   the restriction of $f$ to a leaf $L$ of $\F$  is a homeomorphism onto a leaf $\cL$ of $\F_M$, in the restricted topology on $\cZ$.

 Choose a good covering $\ds \{ \phi_{\alpha} \colon V_{\alpha} \to (-1,1)^n \times (-1,1)^q \mid 1 \leq {\alpha} \leq k\}$ for the foliation $\F$ of $M$, as in \cite{CandelConlon2000}, where $n$ is the leaf dimension of $\F$, and $q$ is the codimension of $\F$ in $M$. Set $T_{\alpha} = \phi_{\alpha}^{-1}(\{0\} \times  (-1,1)^q)$, then the union $\ds T = T_1 \cup \cdots \cup T_k$  is a complete  transversal for  $\F$. We can assume without loss of generality that the closures of the transversals are disjoint.  The Riemannian metric on $TM$ restricts to  a Riemannian metric on each $T_{\alpha}$ and thus defines a path-length metric denoted by $d_{T_{\alpha}}$ on each submanifold $T_{\alpha} \subset M$. Extend the metrics on each $T_{\alpha}$ to a metric $d_T$ on $T$,  by declaring $d_T(u,v) = 1$ if $u \in T_{\alpha}$ and $v \in T_{\beta}$ for $\alpha \ne \beta$.
 
 Recall that  for $(\alpha,\beta)$ admissible,  the overlap of plaques in the charts $V_{\alpha}$ and $V_{\beta}$ defines the holonomy map
 $g_{\alpha, \beta}$. The assumption that $\F$ is a $C^1$-foliation implies that $g_{\alpha, \beta}$  is a $C^1$-map  from an open subset of $T_{\alpha}$ to an open set of $T_{\beta}$.    
For each $u \in \Dom(g_{\alpha, \beta})$, let $D_u(g_{\alpha, \beta})$ denote the   matrix of differentials for $g_{\alpha, \beta}$ at $u \in \Dom(g_{\alpha, \beta})$, with respect to the framing of the tangent spaces to the sections $T_{\alpha}$ induced by the coordinate charts. Let $\|  D_u (g_{\alpha, \beta}) \|$ denote the matrix sup-norm of $D_u (g_{\alpha, \beta}) $  with respect to the Riemannian metric induced on the sections. The assumption that we have a good covering implies that the maps $g_{\alpha, \beta}$ admit continuous $C^1$-extensions, so the norms  $\|  D_u (g_{\alpha, \beta}) \|$ have uniform upper bounds for all admissible pairs $(\alpha, \beta)$  and all $u \in \Dom(g_{\alpha, \beta})$. Define:
\begin{equation}
C_{\F}' = \max \left\{ \|  D_u (g_{\alpha, \beta}) \| \mid    (\alpha, \beta)~ {\rm admissible} ~, ~ u  \in \Dom(g_{\alpha, \beta})  \right\} ~ < ~ \infty
\end{equation}
It follows that the pseudogroup for $\F$ defined by the maps $\{g_{\alpha, \beta}  \mid    (\alpha, \beta)~ {\rm admissible}  \}$  is $C_{\F}'$-Lipschitz.

Recall that $\cT_i \subset \fM$, for $1 \leq i \leq \nu$, are the Cantor transversals to $\fM$ defined by a good covering for $\fM$, as in Definition~\ref{def-fs}.
For each $x \in \cT_i$ there exists $1 \leq \alpha \leq k$ with $f(x) \in V_{\alpha}$, and thus a clopen neighborhood $W(i,x,\alpha) \subset \cT_i$ for which $f(W(i,x,\alpha)) \subset V_{\alpha}$. If $W(i,x,\alpha)$ is sufficiently small, then the plaque projection of the image, $\pi_{\alpha} \colon f(W(i,x,\alpha)) \to T_{\alpha}$, is a homeomorphism onto its image,  
 and so the metric $d_{T_{\alpha}}$ on $T_{\alpha}$ induces   a metric on $W(i,x,\alpha)$. As each $\cT_i$ is compact, we can choose a finite covering $\{ W_{k} \}$ of the union $\cT = \cT_1 \cup \cdots \cT_{\nu}$  where each $W_{k} =  W(i,x,\alpha)$ for appropriate $(i,x,\alpha)$. 
   It may happen that for $x,y \in W_k$ there is an admissible pair $(i,j)$ for the covering of $\fM$ such that $f(h_{i,j}(x))$ and $f(h_{i,j}(y))$ are not contained in the same foliation chart $V_{\ell}$. However, as there are only a finite number of admissible pairs $(i,j)$ for the covering of $\fM$ by foliation charts, we can refine the finite clopen covering $\{ W_{k} \}$ of $\cT$,  so that this condition is satisfied. Then for each $W_k$ and $x \in W_k$ there is an index $i_{\alpha}$ such that $f(x) \in V_{\alpha}$ and $\pi_{\alpha}(f(x)) \in T_{\alpha}$.

We then obtain a metric $d_{\cT}$ on $\cT$ by setting,  
$$d_{\cT}(x,y) =  d_{T_{\alpha}}(\pi_{\alpha}(f(x)), \pi_{\alpha}(f(y))) \quad {\rm if} \quad x,y \in W_k, $$ 
and $d_{\cT}(x,y) = 1$ otherwise.  The metric $d_{\cT}$ induces a metric     on $\fX$,   denoted by $\dX$. We claim there exists $C_{\F} \geq 1$ such that  the action of $\cGX$ on $\fX$ is $C_{\F}$-Lipschitz for   $\dX$ and the generating set $\{ h_{i,j} \mid (i,j) ~ {\rm  admissible} \}$.

Suppose that $x,y \in W_{k}$, then $f(h_{i,j}(x))$ and $f(h_{i,j}(y))$ are   contained in the same foliation chart $V_{\ell}$ by construction. Note that $x$ and $h_{i,j}(x)$ are contained in the same leaf of $\F$ so their images $f(x)$ and $f(h_{i,j}(x))$ are contained in the same leaf of $\F$.  Thus, there is a plaque chain of length at most $\lambda_{f,x}$ between these two points. 
The same holds for the point $y$, so there is a plaque-chain of length $\lambda_{f,y}$ between $f(y)$ and $f(h_{i,j}(y))$. By the compactness of $\cT$,  there is a uniform upper bound $\lambda_f$ for all such pairs.  Thus, by Lemma~\ref{lem-lipalpha} we have the estimate for $x,y \in U_{k}$ with projections $w = \pi(x), w' = \pi(y) \in \fX_i$, and $C_{\F}'' = (C_{\F}')^{\lambda_f}$, 
\begin{equation}\label{eq-Lipschitzlambda}
(C_{\F}'')^{-1} \cdot \dX(w,w') \leq \dX(h_{i,j}(w), h_{i,j}(w')) \leq  C_{\F}'' \cdot \dX(w,w') ~ .
\end{equation}
If $x,y$ do not belong to the same clopen set $W_{k}$, then $\dX(w,w') = 1$ by definition, so their exists   $C_{\F}''' \geq 1$ such that  \eqref{eq-Lipschitzlambda} holds for 
such pairs. Set $C_{\F} = \max \{C_{\F}'', C_{\F}'''\}$, and the claim follows. 
\endproof
 
 \medskip

We next give some properties of Lipschitz pseudogroups and their entropy.  
 The following is an immediate consequence of the definitions.
\begin{lemma}\label{lem-lipalpha}
Suppose that the action of  $\cGX$ on $\fX$ is $C$-Lipschitz  with respect to $\dX$. 
Then for all  $g \in \cGX^*$ with word length $\| g \| \leq \alpha$, and $w,w' \in \Dom(g)$ we have
\begin{equation}\label{eq-Lipschitzalpha}
C^{-\alpha} \cdot \dX(w,w') ~ \leq ~ \dX(g(w), g(w')) ~ \leq ~  C^{\alpha} \cdot \dX(w,w') ~ .
\end{equation}
\end{lemma}

We recall an application of Proposition~2.7 in \cite{GLW1988}, which gives conditions for $h(\cGX,\dX) < \infty$.   

  \begin{prop} \label{prop-entfinite} 
  Let $\fX \subset \mR^q$ be an embedded Cantor set, with metric $\dX$ obtained by the restriction of the standard metric on $\mR^q$. Let $\cGX$ be a finitely generated pseudogroup, with generators $\{h_i \colon U_i \to V_i \mid 1 \leq i \leq k\}$, such that each $h_i$ is the restriction of a $C^1$ diffeomorphism defined on an open neighborhood in $\mR^q$ of the compact set $U_i$.  Then $\cGX$ with the metric $\dX$ is Lipshitz, and the geometric entropy  $h(\cGX, \dX) < \infty$. 
   \end{prop}

  \begin{cor} \label{cor-lipembed}
  Let $\fM$ be a matchbox manifold which embeds as an exceptional minimal set for $C^1$-foliation $\F$ on a compact smooth manifold $M$, as in Proposition~\ref{prop-lipembed}.   Then there is a transverse metric $\dX$ on $\fX$ such that 
 $h(\cGX, \dX) < \infty$. 
   \end{cor}
\proof
Let $\dX$ be the metric on $\fX$ constructed in the proof of Proposition~\ref{prop-lipembed}. Then $\fX$ is covered by disjoint clopen sets for which $\dX$ is the pull-back of the metric on transversals to the foliation $\F$, so by Proposition~\ref{prop-metricentropy2} the entropy for the pseudogroup defined by $\F_M$ restricted to the image of $\fM$ and the entropy for $\cGX$ are either both zero, finite or infinite. Proposition~\ref{prop-entfinite} implies that both entropies are either zero or finite. 
\endproof

Note that by Proposition~\ref{prop-metricentropy1}, the entropy $h(\cGX, \dX)$ is independent of the choice of metric $\dX$ chosen for $\fX$, as long as it defines the topology for $\fX$. Thus by Corollary~\ref{cor-lipembed}  we have:
  \begin{cor} \label{cor-lipembed2}
 Let $\fM$ be a matchbox manifold with pseudogroup $\cGX$ for some regular covering of $\fM$. Suppose there exists a metric $\dX$ on $\fX$ for which 
  $h(\cGX, \dX) = \infty$,   then  $\fM$ is not homeomorphic to an invariant set for any  $C^1$-foliation. 
   \end{cor}

It is well-known that entropy   of a smooth non-singular flow   on  a compact manifold, when restricted to a compact invariant set $\cZ \subset M$, is related to the Hausdorff dimension of $\cZ$, as in \cite{LY1985a,LY1985b}. 

For a Lipschitz pseudogroup $h(\cGX, \dX)$, the box and Hausdorff dimension of $\fX$ with respect to $\dX$ are both well-defined, as in \cite{Edgar1990}, and they do not depend on the Lipschitz equivalence class of the metric $\dX$. While there is no known   direct relation between these dimensions and $h(\cGX, \dX)$, corresponding to the results for flows, there is a finiteness result based on a  concept related to volume doubling for metric spaces (see \cite{Assouad1983,BonkSchramm2000,BuyaloSchroeder2007}).

\begin{defn}\label{def-doubling}
A complete metric space $(\fX, \dX)$ has the \emph{doubling property}, if there exists a constant $C > 1$, such that for every $x \in \fX$, $r > 0$, and integer $n > 0$, the closed ball $B_{\fX}(x,r)$ of radius $r$ about $x$ admits a covering by $C^n$ balls of radius $r/2^n$.
\end{defn}
Note that if $(\fX, \dX)$  has the doubling property, then it has finite box dimension as well.

The proof of  \cite[Proposition~13.2.14]{CandelConlon2000} adapts directly to give:

 \begin{prop} \label{prop-entfinite3} 
If  $h(\cGX, \dX)$ is a compactly generated Lipschitz pseudogroup such that $(\fX, \dX)$ has the doubling property, then $h(\cGX, \dX) < \infty$. 
    \end{prop}

Thus, one approach to constructing matchbox manifolds which cannot embed into a smooth foliation, is to consider examples for which the transversal model space $\fX$ has infinite box dimension, for some metric. This will be discussed further in Section~\ref{sec-nonembedding}.

The   Hausdorff dimension of the transversal  Cantor set to an exceptional minimal set for a $C^1$-foliation is well-studied, especially for foliations of  codimension-one as in     Cantwell and  Conlon \cite{CC1988}, Matsumoto \cite{Matsumoto1988}, Gelfert and Rams  \cite{GelfertRams2009}, and Bi\'{s} and Urbanski \cite{BisUrbanski2008}. 
 Hausdorff dimension is also well-studied for    Cantor sets    defined by contracting Iterated Function Systems (or \emph{IFS}'s),  and the more general   class of self-similar fractals.  For example, see the works  of Rams and his coauthors in \cite{CrovisierRams2006,GelfertRams2009}, and the works of  Rao, Ruan, Wang  and  Xi as in \cite{RRX2006,RRW2012}, are closely related to the study of Lipschitz geometry of foliation minimal sets.
 
 \begin{prob}
 Let $\fM$ be a Lipschitz matchbox manifold, with induced Lipschitz pseudogroup $(\cGX, \dX)$, and suppose  $0 < h(\cGX, \dX) < \infty$. Find properties of the Lipschitz geometry of $\fX$ which must be satisfied if the metric $\dX$ is induced by an embedding of $\fM$ into a $C^r$-foliation, for $r \geq 1$. 
 \end{prob}
  
One can also define finer metric conditions on the action of  a pseudogroup $\cGX$, such as the   Zygmund condition used in \cite{HK1990} which can be used to define ``quasi-conformal'' properties of   homeomorphisms, as in \cite{GardnerSullivan1992,MackayTyson2010,Pansu1989,TukiaVaisala1980,TysonWu2006}.  The study of the Lipschitz properties of Gromov hyperbolic groups acting on their boundaries is a very well-developed subject; see for example   \cite{BuyaloSchroeder2007,KapovichBenakli2001}.

\section{Examples from foliations} \label{sec-foliations}

In this section,  we recall some  examples of    minimal matchbox manifolds which are realized as  exceptional minimal sets in  $C^r$-foliations, for $r \geq 1$. We first consider the   case  for foliations of codimension-one, for which the strongest results have been proven. 
The prototypical example is the well-known  construction by Denjoy:
\begin{thm}[Denjoy \cite{Denjoy1932}]  \label{thm-denjoy}
There exist a $C^1$-diffeomorphism $f$ of the circle $\mS^1$ with no fixed points,  and with a non-empty wandering set $W$ such that the induced action of $f$ on the complement     $\bK = \mS^1 -  W$ gives a minimal   action,  $\vp \colon \mZ \times \bK \to \bK$, called a \emph{Denjoy minimal system}.
\end{thm}

The $C^1$-hypotheses on the diffeomorphism $f$ is far from optimal.   For example,  McDuff \cite{McDuff1981} formulated a set of necessary and sufficient conditions on an embedded Cantor set $\bK \subset \mS^1$ so that it is an invariant set of a $C^{1+\alpha}$-diffeomorphism $f \colon \mS^1 \to \mS^1$, for $0 < \alpha <1$.  
 Other optimal conditions on the derivative of a diffeomorphism $f \colon \mS^1 \to \mS^1$ such that it admits a Cantor minimal set are  discussed in Hu and Sullivan \cite{HuSullivan1997}.

The Denjoy example played a fundamental role in the construction of counter-examples to the Seifert Conjecture, which enabled     Schweitzer in \cite{Schweitzer1974} to construct the first $C^1$-examples of flows on $3$-manifolds without periodic orbits. Schweitzer's construction embedded a suspension of the Denjoy minimal set as an isolated minimal set for a flow contained in a plug embedded in $\mR^3$, and motivated Harrison's construction    \cite{Harrison1988,Harrison1989} of a $C^{2+\alpha}$-flow in $\mR^3$ with an \emph{isolated} minimal limit set homeomorphic to a suspension of the Denjoy set, for $\alpha < 1$.   On the other hand, Knill  constructed in  \cite{Knill1981}  a smooth diffeomorphism in the $2$-dimensional annulus    with a   minimal set homeomorphic to the Denjoy set, so the suspension of this diffeomorphism yields a codimension-$2$ smooth foliation defined by a flow, with a minimal set homeomorphic to the Denjoy minimal set in $\mT^2$. 
Note that   the periodic orbits for the Knill diffeomorphism contain the Denjoy set in its closure, so this example   is not sufficient for constructing smooth counter-examples to the Seifert Conjecture.
The Knill example    illustrates that   the degree of differentiability $r$ for a $C^r$-embedding of a Cantor minimal system may depend on the codimension, as well as the dynamical behavior of the action in open neighborhoods.

In some cases, there are analogs of the above results   for the case of a finitely-generated group acting minimally on a Cantor set. For example, Pixton gave a generalization of the Denjoy construction:
\begin{thm}[Pixton \cite{Pixton1977}] \label{thm-pixton}
Suppose that $0 < \alpha < 1/(n+1)$, then there exist a $C^{1+\alpha}$-action of $\mZ^n$ on the circle $\mS^1$ with no fixed points and with a non-empty wandering set $W$ so that the complement $\bK = \mS^1 -  W$ is a Cantor set which is   minimal for the restricted action.
\end{thm}
   The    Pixton-type  examples have been further studied by Deroin, Kleptsyn and Navas in \cite{DKN2007}, and  Kleptsyn and Navas in \cite{KleptsynNavas2008}.  
   Note that    the suspension of such actions of $\mZ^n$ on $\mS^1$ yield foliations with   exceptional minimal sets, whose leaves are diffeomorphic to $\mR^n$.  
   
Sacksteder proved in    \cite{Sacksteder1965} that  if $\cZ \subset M$ is an exceptional minimal set for a codimension-one $C^2$-foliation $\F_M$ of a compact manifold $M$, then some leaf in $\cZ$ must have an element of  holonomy which is a transverse contraction, and thus cannot be of ``Denjoy type''.  
 A special class of such  examples,  the \emph{Markov minimal sets}, were  studied by Hector \cite{HecHir1981,Hector1983},    
   Cantwell and Conlon  \cite{CC1988}, and Matsumoto \cite{Matsumoto1988}. 
   It remains an open problem to characterize the   embeddings of Cantor minimal systems in $C^r$-foliations of codimension-one, for $r \geq 1$  (see  \cite{Hurder2002}).  
  
There are various constructions of  $C^r$-foliations of codimension $q \geq 2$ with   minimal sets which are matchbox manifold. 
Given a finitely-generated group $\G$ and a $C^r$-action $\vp \colon \G \times N \to N$ on a compact  manifold $N$ of dimension $q$, the suspension of the action (see \cite{CN1985}) yields a $C^r$-foliation of codimension-$q$. In general, it is impossible to determine if such an action $\vp$ has an invariant Cantor set on which the action is minimal, except in very special cases. For example, consider  a lattice subgroup $\G \subset G$ of the rank one connected Lie group $G = SO(q,1)$. The    boundary at infinity for the associated symmetric space $\mH^q = SO(q,1)/O(q)$  is diffeomorphic to  $\mS^q$. If the group $\G$ is a non-uniform lattice,   then the action of $\G$ on its limit set in $\mS^q$ defines a minimal Cantor action, and the suspension of this action is a minimal matchbox manifold embedded in the smooth foliation associated to the action of $\G$ on $\mS^q$.

The   Williams solenoids were introduced in the papers \cite{Williams1967,Williams1974}. Williams proved that for an Axiom A diffeomorphism $f\colon M \to M$ of a compact manifold $M$ with an expanding attractor $\Omega \subset M$, then $\Omega$ admits a stationary presentation, as defined in the next section, and so is homeomorphic to a generalized solenoid. 
  The   unstable manifolds  for $f$ restricted to an open neighborhood  $U$ of $\Omega$  form a   $C^{0,\infty}$-foliation of $U$.  That is, the foliation has $C^0$-pseudogroup maps, with smoothly embedded leaves, and $\Omega$ is the unique minimal set.

\section{Solenoids} \label{sec-solenoids}

 In this section, we describe the constructions of  \emph{weak},   \emph{normal}  and \emph{generalized solenoids}, and   recall   some of their properties.
  We also give a    construction of  metrics on the transverse Cantor sets for which the   holonomy action is by isometries, and hence equicontinuous.   There are many open questions about when such examples can be realized as exceptional minimal sets for $C^r$-foliations.

  A  \emph{presentation}   is a collection $\cP = \{ p_{\ell+1} \colon M_{\ell+1} \to M_{\ell} \mid \ell \geq 0\}$, where each $M_{\ell}$ is a connected compact simplicial complex of dimension $n$, and each  \emph{bonding} map $p_{\ell +1}$  is a proper surjective map of   simplicial complexes with discrete fibers.
For   $\ell \geq 0$ and $x \in M_{\ell}$,  the set $\{p_{\ell +1}^{-1}(x) \} \subset M_{\ell +1}$  is compact and discrete, so   the cardinality $\# \{p_{\ell +1}^{-1}(x) \} < \infty$. It  need not be constant in   $\ell$ or $x$.

    Associated to the presentation $\cP$ is an inverse limit space, called a \emph{generalized solenoid},   
\begin{equation}\label{eq-presentationinvlim}
\cS_{\cP} \equiv \lim_{\longleftarrow} ~ \{ p_{\ell +1} \colon M_{\ell +1} \to M_{\ell}\} ~ \subset \prod_{\ell \geq 0} ~ M_{\ell} ~ .
\end{equation}
 By definition, for a sequence $\{x_{\ell} \in M_{\ell} \mid \ell \geq 0\}$, we have 
\begin{equation}\label{eq-presentationinvlim2}
x = (x_0, x_1, \ldots ) \in \cS_{\cP}   ~ \Longleftrightarrow  ~ p_{\ell}(x_{\ell}) =  x_{\ell-1} ~ {\rm for ~ all} ~ \ell \geq 1 ~. 
\end{equation}
The set $\cS_{\cP}$ is given  the relative  topology, induced from the product topology, so that $\cS_{\cP}$ is itself compact and connected.

For example, if    $M_{\ell} = \mS^1$ for each $\ell \geq 0$, and the map $p_{\ell}$ is a proper covering map of degree $m_{\ell} > 1$ for $\ell \geq 1$, then $\cS_{\cP}$ is an example of a  {classic solenoid},  discovered independently by    van~Dantzig \cite{vanDantzig1930} and   Vietoris   \cite{Vietoris1927}.

   We say the presentation $\cP$ is \emph{stationary} if  $M_{\ell} = M_0$ for all $\ell \geq 0$, and the bonding maps $p_{\ell} =p_1$ for all $\ell \geq 1$. A solenoid $\cS_{\cP}$ obtained from a stationary presentation $\cP$ has a self-map $\sigma$ defined by the shift, 
$\ds  \sigma(x_0, x_1, \ldots ) = (x_1, x_2, \ldots )$. The map $\sigma$  can be considered as a type of expanding map on  $\cS_{\cP}$, though in fact it may   be expanding only in some directions, as discussed in Section~3 of \cite{BHS2006}. By the work of Mouron \cite{Mouron2009,Mouron2011}, these are the only examples of $1$-dimensional solenoids   with an expanding map. The case for expanding maps of generalized $1$-dimensional solenoids is much richer, as described in the work of Williams \cite{Williams1967,Williams1970}, which classifies the stationary inverse limits defined by   expanding maps of branched $1$-manifolds.

If     $M_{\ell}$ is a compact manifold without boundary for each $\ell \geq 0$,  and  the map $p_{\ell}$ is a proper covering map of degree $m_{\ell} > 1$  for   $\ell \geq 1$,   then $\cS_{\cP}$  is said to be a \emph{weak solenoid}. This generalization of $1$-dimensional solenoids was originally    considered  in the papers by McCord \cite{McCord1965} and Schori  \cite{Schori1966}. In particular, McCord showed in \cite{McCord1965} that   $\cS_{\cP}$ has a local product structure, hence

\begin{prop}\label{prop-solenoidsMM}
Let   $\cS_{\cP}$ be a weak solenoid, whose    base space $M_0$   is a compact manifold of dimension $n \geq 1$. Then   $\cS_{\cP}$ is  a minimal matchbox manifold of dimension $n$.
\end{prop}

Associated to a presentation $\cP$ of compact manifolds is a sequence of proper surjective maps 
$$q_{\ell} = p_{1} \circ \cdots \circ p_{\ell -1} \circ p_{\ell} \colon M_{\ell} \to M_0 ~ .$$
For each $\ell > 1$, projection onto the $\ell$-th factor in the product $\ds \prod_{\ell \geq 0} ~ M_{\ell}$ in \eqref{eq-presentationinvlim} yields a 
  fibration map denoted by $\Pi_{\ell} \colon \cS_{\cP}  \to M_{\ell}$, for which 
 $\Pi_0 = \Pi_{\ell} \circ q_{\ell} \colon \cS_{\cP} \to M_0$. 
A choice of a basepoint $x \in \cS_{\cP}$ gives basepoints 
$x_{\ell} = \Pi_{\ell}(x) \in M_{\ell}$, and we define $\cH^x_{\ell} = \pi_1(M_{\ell}, x_{\ell})$. Let  $\fX_x = \Pi_0^{-1}(x)$ denote the fiber of $x$, which is Cantor set by the assumption on the cardinality of the fibers of each map $p_{\ell}$. 

 A presentation $\cP$  is said to be \emph{normal} if, given a basepoint $x \in \cS_{\cP}$,   for each $\ell \geq 1$  the image  subgroup of the map 
$\ds (q_{\ell} )_{\#} \colon \cH^x_{\ell} \longrightarrow \cH^x_{0}$ 
is    a normal subgroup. 
Then each quotient  $G^x_{\ell} = \cH^x_{0}/\cH^x_{\ell}$ is      finite group, and there are surjections $G^x_{\ell +1} \to G^x_{\ell}$. The fiber $\fX_x$   is then naturally identified with the \emph{Cantor group} defined by the inverse limit,  
\begin{equation}\label{eq-Galoisfiber}
G^x_{\infty} = \lim_{\longleftarrow} ~ \{ p_{\ell +1} \colon G^x_{\ell +1} \to G^x_{\ell }\} ~ \subset \prod_{\ell \geq 0} ~ G^x_{\ell} ~ .
\end{equation}
The   fundamental group $\cH^x_0$ acts on the fiber $G^x_{\infty}$ via     the coordinate-wise multiplication on the product in \eqref{eq-Galoisfiber}.  In the case of the Vietoris solenoid, where each map $p_{\ell} \colon \mS^1 \to \mS^1$ is a double cover, the fiber $G^x_{\infty}$ is the dyadic group. 
More generally, 
a   solenoid $\cS_{\cP}$ is said to be a \emph{normal} (or \emph{McCord}) \emph{solenoid} if the tower of coverings in the presentation is normal, and thus the fiber over $x \in M_0$ of the map $\cS_{\cP} \to M_0$ is the Cantor group $G^x_{\infty}$.

\begin{lemma}\label{lem-denseaction}
Let $\cP$ be a presentation of a weak solenoid $\cS_{\cP}$, choose a basepoint $x \in \cS_{\cP}$ and set $\fX_x = \Pi_0^{-1}(x)$, and recall that $\cH_0^x = \pi_1(M_0,x_0)$.   Then   the left action of $\cH_0^x$ on $\fX_x$ is minimal. 
\end{lemma}  
 \proof
 The left action of $\cH_0^x$ on each quotient space $X_{\ell} = \cH^x_{0}/\cH^x_{\ell}$ is transitive, so the orbits are dense in the product topology for $\fX_x$.
  \endproof

Let $\wtM_0$ denote the universal covering of the compact manifold $M_0$. Associated to the left action of $\cH_0^x$ on $\fX_x$ is a suspension minimal matchbox manifold 
\begin{equation}\label{eq-suspensionfols}
\fM = \wtM_0 \times \fX_x / (y_0 \cdot g^{-1}, x) \sim (y_0 , g \cdot x) \quad {\rm for }~ y_0 \in \wtM_0 , ~ g \in \cH_0^x ~.
\end{equation}
Given  coverings $\pi' \colon M' \to M$ and $\pi'' \colon M'' \to M$, such that    the subgroups 
$$\ds \pi_{\#}'(\pi_1(M',x')) = \pi_{\#}''(\pi_1(M'', x'')) \subset \pi_1(M,x) ,$$ 
then there is a natural homeomorphism of coverings $M' \cong M''$ which is defined using the path lifting property. From this, it easily follows (see   \cite{ClarkHurder2013}) that:
\begin{prop}\label{prop-weaksuspensions}
Let   $\cS_{\cP}$ be a weak solenoid with base space $M_0$ where $M_0$ is a compact manifold of dimension $n \geq 1$. Then there is a foliated homeomorphism $\cS_{\cP} \cong \fM$.
\end{prop}
 
 \begin{cor}\label{cor-weaksuspensions}
The homeomorphism type of a weak solenoid    $\cS_{\cP}$ is completely determined by the base manifold $M_0$ and the descending chain of subgroups   
\begin{equation}
 \cH^x_{0} ~ \supset  ~  \cH^x_{1} ~ \supset  ~  \cH^x_{2} ~ \supset  ~  \cH^x_{3} ~ \supset  ~ \cdots 
\end{equation}
\end{cor}

 Note      the   intersection $\ds \cH^x_{\infty} \equiv \cap_{\ell \geq 1} ~ \cH^x_{2}$ is the fundamental group of the typical leaf of $\cP$. If this intersection group is trivial , then all leaves of the foliation $\F$ for  
$\fM \cong \cS_{\cP}$ are isometric to the universal covering of the base manifold $M_0$.

 The presentation $\cP$ of an inverse limit    $\cS_{\cP}$ can be used to construct a ``natural'' metric on the space, and which is well-adapted to Lipschitz maps between such spaces. This has been studied in detail in the works by Miyata and Watanabe \cite{MiyataWatanabe2002,MiyataWatanabe2003a,MiyataWatanabe2003b,MiyataWatanabe2003c,Miyata2009}. In the case of weak solenoids, this construction of natural metrics adapted to the resolution takes on a simplified form.

Let  $\cS_{\cP}$ be a weak solenoid, with notations as above. Then each quotient $X_{\ell} = \cH^x_{0}/\cH^x_{\ell}$ is a finite set with a transitive left action of the fundamental group $\cH^x_{0}$. Let  $d_{\ell}$ denote the discrete metric on $X_{\ell}$, where $d_{\ell}(x,y) =1$ unless $x=y$, for $x,y \in X_{\ell}$.   Observe that the left action of $\cH^x_{0}$   acts by isometries for the metric $d_{\ell}$.  Choose a  positive series $\{a_{\ell} \mid a_{\ell} > 0\}$ with total sum $1$, then define a metric on $\fX_x$ by setting, for $u,v \in \fX_x$ so 
$u = (x_0, u_1, u_2, \ldots)$ and $v = (x_0, v_1, v_2, \ldots)$, 
\begin{equation}\label{eq-canonicalmetric}
\dX(u,v) = a_1 d_1(u_1, v_1) + a_2 d_2(u_2 , v_2) + \cdots 
\end{equation}
Then $\dX$  is invariant under the action of $\cH^x_0$, so the holonomy   for  the fibration $\Pi_0 \colon \cS_{\cP} \to M_0$ acts by isometries for this metric on $\fX_x$.

It may happen that we have two presentations $\cP$ and $\cP'$ over the same base manifold $M_0$ such that their inverse limits are homeomorphic as fibrations,  $h \colon \cS_{\cP} \cong \cS_{\cP'}$. However, the map $h$ need not be Lipschitz on fibers for the metrics associated to the presentations as above, as will be seen in the examples in Section~\ref{sec-classification}.

The normal solenoids have a nice characterization among the matchbox manifolds. A continuum $\Omega$ is \emph{homogeneous} if its group of homeomorphisms is  transitive. That is, given any two points $x,y \in \fM$, there is a   homeomorphism $\ds h \colon \fM \to \fM$ such that $h(x) = y$.  It was shown in   \cite{ClarkHurder2013} that:
\begin{thm}\label{thm-homogeneous}
Let $\fM$ be a  homogeneous matchbox manifold. Then   $\fM$ is homeomorphic to a normal solenoid  $\cS_{\cP}$ as foliated spaces.
\end{thm}

The normal solenoids are the analogs in codimension-zero foliation theory, for  the transversely parallelizable (TP) equicontinuous foliations   in a topological version of Molino's Theory for smooth foliations of   manifolds \cite{ALMG2013}.
Note that all leaves in a normal solenoid are homeomorphic, as the spaces are homogeneous.  
In the case of weak solenoids, the leaves of $\F$ need not be homeomorphic, and the works \cite{CFL2010,DDMN2010} give examples where the leaves of $\F$  have differing numbers of ends. There is no analog of this behavior in the context of smooth Riemannian foliations   on manifolds.

Now consider a matchbox manifold $\fM$ of dimension $n$, but whose associated pseudogroup $\cGX$ is not  equicontinuous.   This type of matchbox manifold arises  in the study of   the tiling spaces associated to aperiodic tilings of $\mR^n$ with finite local complexity, and also as foliation minimal sets. For example,   the Hirsch examples in  \cite{Hirsch1975} (see also \cite{BHS2006}) yield  real analytic foliations of codimension-one  with    exceptional minimal sets and expansive holonomy pseudogroups. 
 Also, the exceptional minimal sets for the Denjoy and Pixton examples discussed in Section~\ref{sec-foliations} have the property that all of their leaves are diffeomorphic to $\mR^n$, and so they are without leafwise holonomy, but the global holonomy pseudogroup $\cGX$ associated to them is not equicontinuous.  
  It follows from the following   result that  each of their minimal sets    admits a presentation of the form \eqref{eq-presentationinvlim}.

\begin{thm}[\cite{CHL2013b}]  \label{thm-shapemm}
Let $\fM$ be a minimal  matchbox manifold without germinal  holonomy. Then there exists a presentation $\cP$ by  simplicial maps   between compact branched manifolds,   such that $\fM$ is homeomorphic to $\cS_{\cP}$ as foliated spaces.  
\end{thm}

 \begin{cor}  \label{cor-shapemm}
Let $\fM$ be an exceptional  minimal  set for a $C^1$-foliation $\F$ of a compact manifold $M$.  If all leaves of $\F | \fM$ are simply connected, then  there is a    homeomorphism of $\fM$ with the inverse  limit space $\cS_{\cP}$ defined by   a   presentation $\cP$, given by  simplicial maps   between compact branched manifolds.
\end{cor}
 
In the case of the Denjoy and Pixton examples given in Theorems~\ref{thm-denjoy} and \ref{thm-pixton}, the   geometry of their construction implies that the presentation $\cP$ one obtains is stationary. 

\begin{prob}\label{prob-stationary}
Let $\fM$ be an exceptional  minimal  set for a $C^r$-foliation $\F$ of a compact manifold $M$, where $r \geq 1$, and assume that $\fM$ is without holonomy. Find conditions on the holonomy pseudogroup $\cGX$ for $\F$ which are sufficient to imply that $\fM$ admits a stationary presentation.
\end{prob}
One approach to this problem, is to ask if the existence of approximations to the foliation $\F$ on $\fM$ by the compact branched manifolds $M_{\ell} = M_0$ in a stationary presentation $\cP$, implies some form of ``finiteness'' for the holonomy maps of the pseudogroup $\cGX$.  Such finiteness  conditions may be derived, for example,  from the induced action of the   shift map $\sigma$ on  the   tower of maps in the presentation. Then one would try to ``fill in'' the approximations with a foliation on an open neighborhood.  Such a result would be reminiscent of the approach to showing the vanishing of the Godbillon-Vey class by   Duminy and Sergiescu  in \cite{DS1981}. 

 Theorem~\ref{thm-shapemm} is a generalization of a celebrated result by Anderson and Putnam in \cite{AP1998} for tiling spaces.  Given a   repetitive, aperiodic tiling of the Euclidean space $\mR^n$ with finite local complexity, the associated tiling space $\Omega$ is defined as the   closure of the set of translations by $\mR^n$ of the given tiling,   in an appropriate    topology  on the  space of tilings on $\mR^n$. The space $\Omega$ is a matchbox manifold in our sense, whose leaves are defined by a free action of $\mR^n$ on $\Omega$ (for example, see \cite{PFS2009,SW2003,Sadun2008}.)
A remarkable result in the theory of tilings of  $\mR^n$ is  that  the  tiling space $\Omega$  admits a presentation  as the inverse limit of a tower of branched flat manifolds \cite{AP1998,Sadun2003, Sadun2008}, where the   branched manifolds are the union of finite collections of tiles.

Other generalizations of the Anderson-Putnam theorem   have been given. For example,  the work of Benedetti and Gambaudo   in \cite{BG2003} discusses constructing towers for special classes of matchbox manifolds \emph{with possibly non-trivial but finite holonomy}, where the leaves are defined by a locally-free action of a connected Lie group $G$.  Their work suggests what appears to be a difficult problem:

\begin{prob}\label{prob-torsion}
Let $\fM$ be a minimal matchbox manifold with leaves having non-trivial holonomy. Show that $\fM$ is homeomorphic to an inverse limit $\cS_{\cP}$ for some modified notion of  presentations by branched manifolds, which takes into account   the leafwise holonomy groups.   
\end{prob}
Note that a solution to this problem would yield a presentation for  an exceptional minimal set in a $C^2$-foliation of codimension-one, which  by the results of Sacksteder in \cite{Sacksteder1965}  always have leaves with holonomy. The existence of such a presentation   would provide an alternate approach to the celebrated result of Duminy on the ends of leaves in exceptional minimal sets \cite{CantwellConlon2002}.

Theorem~5.8 in the paper \cite{LR2013}  states a solution to   Problem~\ref{prob-torsion}, though it seems that the claimed result conflicts with the results of  \cite{BG2003} for a model of generalized tiling spaces defined by $G$-actions with non-trivial holonomy. Also, the results of Section~6 of the same paper conflict with other established results concerning weak solenoids.

  \begin{prob}\label{prob-weaksolenoids}
Given a weak solenoid $\cS_{\cP}$ with presentation $\cP$ and associated transverse metric   given by \eqref{eq-canonicalmetric}, does there exists a Lipschitz embedding of $\cS_{\cP}$ as an exceptional minimal set for a $C^r$-foliation of a smooth manifold $M$?
\end{prob}
The problem is of interest whether $M$ is assumed compact, or open without boundary, and for any $r \geq 1$. All known results are for the case where the base $M_0 = \mT^n$ is a torus, for $n \geq 1$. The examples of type DE (\emph{Derived from Expanders}) described by Smale in \cite[p. 788]{Smale1967}, constructs an embedding of the dyadic solenoid over $\mS^1$ which is realised as a basic set for a smooth diffeomorphism and is  an attractor. More general realizations of $1$-dimensional solenoids   as minimal sets for  smooth flows were  constructed for  flows in the works by Gambaudo and Tresser \cite{GST1994}, Gambaudo, Sullivan and Tresser \cite{GST1994}, and Markus and Meyer \cite{MM1980}.  
The case when the base manifold $M_0 = \mT^n$ for $n \geq 2$ was studied by the author with Clark in \cite{ClarkHurder2011}, where it was shown that for every presentation $\cP$ there exists a refinement $\cP'$ which can be realized in a $C^r$-foliation. That is, every topological type can be realized, though the metric induced on the inverse limit depends on the presentation $\cP$.

All of the known examples of weak solenoids which embed as exceptional minimal sets for $C^2$-foliations have abelian fundamental group $\cH_x$ and so are consequently normal solenoids.  It seems plausible, based on the proofs in \cite{ClarkHurder2011}, to conjecture that if a weak solenoid admits   an embedding  in a $C^2$-foliation, then  it must be a normal solenoid with nilpotent  covering groups. It also seems possible that an even  stronger conclusion holds, that the covering groups for such a smoothly embedded solenoid must be abelian.

\section{Fusion of Cantor minimal systems} \label{sec-fusion}

There is a well-known method, called  \emph{tubularization}, of amalgamating the holonomy pseudogroups of   two     foliations $\F_1, \F_2$  of codimension-one with  the same leaf dimension. We recall this method briefly,    then introduce the analogue of  this  technique for minimal matchbox manifolds, to obtain the \emph{fusion} of their holonomy pseudogroups. 
 
Assume there are given two foliations say $\F_1$ and $\F_2$, on manifolds $M_1$ and $M_2$, of with leaf dimension $n$ and codimension-one. We assume that  their normal bundles are oriented,   and there are given  smooth embeddings  $\eta_i \colon \mS^1 \to M_i$ which are transverse to $\F_i$ for $i=1,2$. 
For $\e > 0$   small, let $\cE(\eta_i, \e) \subset M_i$ be the  closed $\e$-disk neighborhood of the image of the map $\eta_i$, where we assume $\e> 0$ is chosen so that $\cE(\eta_i, \e)$  is an embedded submanifold with boundary diffeomorphic to $\mT^2$. Then the restriction of $\F_i$ to $\cE(\eta_i, \e)$ is a foliation whose leaves   are closed $2$-disks, and which are parametrized by   $\mS^1$ via the transversal $\eta_i$. 

 The choice of a diffeomorphism $\vp \colon \mS^1 \to \mS^1$ extends to give a foliated map   $\whvarp \colon \cE(\eta_1, \e) \to \cE(\eta_2, \e)$, which we use to identify   the boundaries $\partial \cE(\eta_1, \e)$ and $\partial \cE(\eta_2, \e)$.  Denote the resulting surgered manifold   by $M = M_1 \#_{\vp} M_2$. 
Then $M$  has a foliation of codimension-one,   whose foliation pseudogroup is the amalgamation, or ``pseudogroup free product'',  of the pseudogroups for $\F_1$ and $\F_2$. This very useful construction has  many   applications \cite{CN1985,CandelConlon2000,HecHir1981,Lawson1977}. 

For foliations with codimension $q > 1$, the tubularization method is not so commonly used, as the existence of a compact manifold $N$ and embeddings $\eta_i \colon N \to M_i$ transverse to the given foliations is a highly   exceptional condition to assume. The tubularization method is   often replaced with the method of \emph{spinnable structures} of Tamura \cite{Tamura1972}, or the  \emph{open book} method as in \cite{Lawson1971,Winkelnkemper1973}.

Next,  we define the analog of tubularization for Cantor pseudogroups. We first describe this construction for group actions. Assume there are given    actions $\vp_i \colon \G_i \times \bK_i \to \bK_i$ for $i=1,2$,   of   finitely generated groups $\G_i$ on   Cantor sets $\bK_i$. 
Choose  clopen subsets $V_i \subset \bK_i$ and a homeomorphism $h \colon V_1 \to V_2$. Define the Cantor set $\bK = \bK_1 \#_h \bK_2$ obtained from the disjoint union $\bK_1 \cup \bK_2$ by identifying the clopen subsets $V_1$ and $V_2$ using the map $h$. 

The action on $\bK$ of  $\gamma \in \G_1$   is  via $\vp_1(\gamma)$ on $\bK_1$, and acts as the identity on the complement $\bK_2 - V_2$. Analogously, the action of $\vp_2$ extends to an action of the elements of $\G_2$ on $\bK$. This produces an action $\vp$ of the free product $\G_1 * \G_2$ on $\bK$. 
Note that if  $V_1 = \bK_1$ and $V_2 = \bK_2$ then this process is just combining the generators of $\vp_1(\G_1)$ with the conjugates by $h$ of the generators of $\vp_2(\G_2)$.

If each of the actions $\vp_i$ is minimal, then the action of $\vp$ on $\bK$ is also minimal.

In the case where $\cG_{\bK_1}$ is a pseudogroup acting on $\bK_1$ and $\cG_{\bK_2}$ is a pseudogroup acting on $\bK_2$, then the amalgamation of their actions over a homeomorphism $h \colon V_1 \to V_2$ is actually simpler, as there is no need to extend the domains of the local actions.

If the action $\vp_i$ is realized as the holonomy of a suspension matchbox manifold $\fM_i$ as in \eqref{eq-suspensionfols}, then the action of $\vp$ is realized as the holonomy of a surgered matchbox manifold $\fM = \fM_1 \#_h \fM_2$ constructed analogously to the method described above for codimension-one foliations.      This construction is analogous to the construction of a new graph matchbox manifold,  from two given graph matchbox manifolds, which was introduced  by Lukina in \cite{Lukina2012} as part of her study of the dynamics of examples obtained by  the Ghys-Kenyon construction. Lukina called this process ``fusion'', and we adopt the same terminology for the process described here.  
 
 \begin{defn}\label{def-fusion}
 Let $\fM_i$ be minimal matchbox manifolds with transversals $\fX_i$ for $i =1,2$. Choose clopen subsets $V_i \subset \fX_i$ and a homeomorphism $h \colon V_1 \to V_2$. 
 Then the minimal matchbox manifold $\fM = \fM_1 \#_h \fM_2$ is said to be the \emph{fusion} of $\fM_1$ with $\fM_2$ over $h$. 
 \end{defn}

The concept of fusion for matchbox manifolds illustrates some of their fundamental differences with smooth foliations.  A clopen transversal for a  smooth foliation must be a compact submanifold without boundary, which does not always exist, while   the above fusion construction  can always be defined, along with many variations of it.
 Here is an interesting basic question:
\begin{prob}\label{prob-fusiondynamics}
How are the dynamical properties of a fusion $\fM = \fM_1 \#_h \fM_2$ related to the   dynamical properties of the factors $\fM_1$ and $\fM_2$? In particular, describe the geometric structure of the leaves in $\fM$, in terms of the structure of the leaves of the factors $\fM_1$ and $\fM_2$ and the fusion map $h \colon V_1 \to V_2$ between transversals. Show that the theory of  hierarchies for the leaves of graph matchbox manifolds in Lukina   \cite{Lukina2012}  also apply for fusion in the context of matchbox manifolds.
\end{prob}

\section{Non-embeddable matchbox manifolds} \label{sec-nonembedding}

In this section, we   construct   examples of   Lipschitz pseudogroups $(\cGX, \dX)$  which cannot arise from an embedding of a matchbox manifold into a $C^1$-foliation. 
  All of the pseudogroups constructed   can be realized as the holonomy of a   matchbox manifold $\fM$, using the suspension construction described in \cite{LRL2013}.
Thus, the resulting matchbox manifolds  $\fM$ do not embed as closed invariant sets for any $C^1$-foliation. There are many variations on the constructions, which shows that  there is a wide variety of non-embeddable   matchbox manifolds. 

 The idea of the construction is to produce a Lipschitz pseudogroup $\cGX$  with infinite entropy, $h(\cGX, \dX) = \infty$, so that by  Corollary~\ref{cor-lipembed} the associated suspension matchbox manifold is not homeomorphic to an exceptional minimal set.     Achieving infinite entropy with Lipschitz generators for $\cGX$  requires that the space $(\fX, \dX)$ have infinite Hausdorff dimension. The first step then, is the construction of the model for the metric Cantor set $(\fX, \dX)$, which is based on   the construction of   \emph{graph matchbox manifolds}, as introduced by Ghys in \cite{Ghys1999}, and studied   in 
 \cite{Blanc2003,LR2011,LRL2013,Lukina2014}.

 Let $\cT$ be an infinite connected tree with bounded valence. The example that we consider here is the Cayley graph $T_n$ for the free group on $n$-generators, $\mF_n = \mF * \cdots * \mF$, for $n \geq 2$. Choose a basepoint $e \in \cT$.   Each edge of  $\cT$ is homeomorphic to $[0,1]$   so inherits a metric from $\mR$. Then give $\cT$ the path length metric, and let $B_{\cT}(x,n) \subset \cT$ denote the closed ball of radius $n$ centered at $x \in \cT$. Thus, if $x$ is a vertex of the tree, then $B_{\cT}(x,n)$ is a connected subtree of $\cT$.
 
We say that a subtree $T \subset \cT$ has a \emph{dead end}, if there is a vertex $x \in T$ which is contained in a unique edge.  Let $X$ be the set of all connected   subtrees of $\cT$ which have no dead ends, and such that $e \in T$.   Define the metric $d_X$    on $X$ by declaring that, for $T, T' \in X$ then 
 $$d_X(T, T') \leq 2^{-n} ~ \Longleftrightarrow ~ B_{\cT}(x,n) \cap T = B_{\cT}(x,n) \cap T'. $$
Let $\fX$ denote   the closure of $X$ in this metric, then $\fX$ is a totally disconnected space. A point $z \in \fX$ is then a subtree of $\cT$ which contains the basepoint $e$.

In the case   where $\cT_n$ is the Cayley graph of $\mF_n$, we denote the closure of the space of subtrees of $\cT_n$ as above by $\fX_n$.  The ``no dead end'' assumption on the subtrees implies that $\fX_n$ has no isolated points, hence is a Cantor set. Let $d_{\fX_n}$ denote the induced metric on $\fX_n$. Then  we have:
 
\begin{thm}[Lukina \cite{Lukina2014}]\label{thm-infdim}
For    $n \geq 2$,   the metric space $(\fX_n , d_{\fX_n})$ has infinite Hausdorff dimension.
\end{thm}

The translation action of $\mF_n$ on $\cT_n$ defines a pseudogroup $\cG_{\fX_n}$ acting on $\fX_n$, where a word $\gamma \in \mF^n$ acts on the pointed subtree $(T,e)$ if $\gamma \cdot e \in T$, so that $(\gamma^{-1} \cdot T, e) \in \fX_n$. This action is discussed further  in  \cite{LRL2013,Lukina2014}. In particular, the action is Lipschitz for the metric $d_{\fX_n}$, with   constant $C = 2$.

Lukina shows in \cite{Lukina2012} that there exists a dense orbit for this action, so the pseudogroup is transitive. However, the periodic orbits for the action of $\cG_{\fX_n}$ are dense, so the action is not minimal. 

 The proof of Theorem~\ref{thm-infdim} in \cite{Lukina2014} essentially   shows the following, with the details given in   \cite{HL2014}:
 
 \begin{thm}\label{thm-infent}
For    $n \geq 2$,  $h(\cG_{\fX_n}, d_{\fX_n}) = \infty$ for the  metric space $(\fX_n , d_{\fX_n})$.
\end{thm}

The suspension construction for pseudogroups given in \cite{LRL2013} constructs a $2$-dimensional matchbox manifold $\fM_n$ whose holonomy pseudogroup is $\cG_n$. Thus, combining  Theorem~\ref{thm-infent} with   Corollary~\ref{cor-lipembed} and Propositions~\ref{prop-metricentropy1} and \ref{prop-metricentropy2}, we have the consequence:

 \begin{thm}\label{thm-noembedMn}
For $n \geq 2$, the transitive Lipschitz matchbox manifold $\fM_n$ is not homeomorphic to an invariant subset of any $C^1$-foliation $\F_M$ of a   manifold $M$. 
\end{thm}

Now consider a minimal Cantor action $\vp_2  \colon \fX_2 \to \fX_2$ for some Cantor set $\fX_2$.  For example, let $\vp_2$ be a Denjoy type homeomorphism. Then there is a homeomorphism $h \colon \fX_n \to \fX_2$ and we can form the fusion of the action of $\cG_{\fX_n}$ with that of $\vp_2$. 
That is, we adjoin the action of $\whvarp_2 \equiv h^{-1} \circ  \whvarp_2 \circ h$ to the action of $\cG_{\fX_n}$ on $\fX_n$ to obtain a minimal action of the fusion  pseudogroup,   denoted by $\widehat{\cG}_{\fX_n}$. Let $\widehat{\fM}_n$ denote the suspension matchbox manifold obtained from $\widehat{\cG}_{\fX_n}$.

The action of $\whvarp_2$ is not assumed to be Lipschitz, but we have in any case:

  \begin{thm}\label{thm-noembedMn2}
For $n \geq 2$, the minimal   matchbox manifold $\widehat{\fM}_n$ is not homeomorphic to an invariant subset of any $C^1$-foliation $\F_M$ of a   manifold $M$. 
\end{thm}
\proof
Suppose that $\widehat{\fM}_n$ is   homeomorphic to an invariant subset $\cZ \subset M$ of a $C^1$-foliation $\F_M$ on $M$. Then $\cZ$ must be a saturated subset, and every leaf is dense as this is true for $\widehat{\fM}_n$. Moreover, the transversals to $\widehat{\fM}_n$ are Cantor sets, so $\cZ$ must be an exceptional minimal set for $\F_M$.  Then by Proposition~\ref{prop-lipembed}, the embedding induces a metric $d_{\fX_n}'$ on $\fX_n$ such that $\widehat{\cG}_{\fX_n}$ is a Lipschitz pseudogroup for this metric.
By construction, $\widehat{\cG}_{\fX_n}$ contains   $\cG_{\fX_n}$ as a sub-pseudogroup, and so 
$$h(\widehat{\cG}_{\fX_n}, d_{\fX_n}') ~ \geq ~  h(\cG_{\fX_n}, d_{\fX_n}') ~ = ~  h(\cG_{\fX_n}, d_{\fX_n}) ~ = ~\infty$$
where we use Proposition~\ref{prop-metricentropy1}. But this contradicts  Corollary~\ref{cor-lipembed}.
\endproof

  These two examples suggests the following:   
  \begin{prob}
 Show that there is no metric  $d_{\fX_n}''$ on $\fX_n$ for which the action of $\widehat{\cG}_{\fX_n}$ is Lipschitz.
  \end{prob}
 It seems very likely that this has a positive solution, that no such metric can exists, though the proof of this fact may require some new insights or techniques.  
   
   \medskip
   
We conclude this section with another remark, and a question. Recall that  Problem~\ref{problem8} asks for  obstructions to the existence of an embedding $\iota \colon \fM \to M$ of a Lipschitz matchbox manifold as an exceptional minimal set for a $C^1$-foliation $\F$ on $M$.  Such an embedding implies in particular that the transverse Cantor set $\fX$ admits a Lipschitz embedding into the Euclidean space $\mR^q$.  The question of when a metric space admits a Lipschitz  embedding in $\mR^q$  dates from the 1928 paper \cite{Bouligand1928}, and is certainly well-studied. For example, the  doubling property in Definition~\ref{def-doubling} of Assouad \cite{Assouad1983}, and the weakening of this condition by Olson and Robinson \cite{OlsonRobinson2010}, prove embedding criteria for metrics.  These are types of  ``asymptotic small-scale homogeneity'' properties of the metric $\dX$, which  suggests   an alternate approach to the Lipschitz embedding problem for minimal pseudogroups. 

\begin{prob} \label{prob-doubling}
Let $\fX$ be a Cantor space with metric $\dX$. Let $\cGX$ be a compactly-generated     pseudogroup acting minimally on $\fX$, and which is Lipschitz with respect to $\dX$. If the metric $\dX$ satisfies some version of the doubling condition, so that $(\fX, \dX)$ admits a Lipschitz embedding into some $\mR^q$, does there also exists an  embedding such that    $\cGX$ is obtained by the restriction of some $C^1$-pseudogroup acting on an open neighborhood of the embedded Cantor set?
\end{prob}
 
 For a Cantor set $\fX$ with an ultrametric $\dX$, the Lipschitz embedding problem for $(\fX, \dX)$ has been solved for various special cases.  The work of  Julien and  Savinien in \cite{JS2011} estimates the Hausdorff dimension   for a self-similar Cantor set with an ultrametric, and they derive    estimates for  its Lipschitz embedding dimension. The embedding properties of ultrametrics on Cantor sets which are the boundary of a hyperbolic group  are discussed  by 
Buyalo and Schroeder in \cite[Chapter 8]{BuyaloSchroeder2007}. In both of these cases, it seems likely that the answer to Problem~\ref{prob-doubling} is positive. In general, one expects the solution to be more complicated, as is almost always the case with Cantor sets.

 Finally, recall  that every Cantor set embeds homeomorphically to a Cantor set in $\mR^2$, and any two such are homeomorphic by a homeomorphism of $\mR^2$ restricted to the set. This classical fact, due to Brouwer, is proved in detail by Moise in Chapter 12 of   \cite{Moise1977}. It has been used to construct topological embeddings of solenoids in codimension-two   foliations, as in the work of Clark and Fokkink \cite{ClarkFokkink2004}.

 On the other hand, the tameness property of Cantor sets in $\mR^2$ does not hold for all  Cantor sets embedded in $\mR^3$. The \emph{Antoine's Necklace}   is the   classical example of this, as discussed in Chapter 18 of \cite{Moise1977}, and in Section~4.6 of \cite{HockingYoung1988}. 
   It seems natural to ask the naive question: 
 
 \begin{prob} \label{prob-antoine}
Let $\fA$  denote the Antoine Cantor set embedded in $\mR^3$, with the   metric $\dA$ on $\fA$ induced by the restriction of the Euclidean metric. Does there is some exceptional minimal set for a $C^1$-foliation of codimension three, whose transverse model space is Lipschitz equivalent to $(\fA, \dA)$?
\end{prob}

\section{Classification of Lipschitz solenoids} \label{sec-classification}

In this section, we define  \emph{Morita equivalence} and  \emph{Lipschitz equivalence} of minimal   pseudogroups, and consider the problem of Lipschitz classification for the special case of normal solenoids. While the condition of Morita equivalence is well-known and studied, Lipschitz equivalence seems less commonly studied, except possibly for group and semi-group actions on their boundaries.

Let $\cGX$ be a minimal pseudogroup acting on a Cantor space $\fX$, and let $V \subset \fX$ be a clopen subset. The induced pseudogroup $\cGX | V$ is defined as the subcollection of all maps in $\cGX$ with domain  and range in $V$.   The following is then the adaptation of the notion of Morita equivalence of groupoids, as  in Haefliger \cite{Haefliger1984},  to the context of minimal Cantor actions.

\begin{defn}
Let $\cGX$ be a minimal pseudogroup action on the Cantor set $\fX$  via Lipschitz homeomorphisms   with respect to the metric $\dX$. 
Likewise, let   $\cGY$ be a minimal pseudogroup action on the Cantor set $\fY$    via Lipschitz homeomorphisms with respect to the  metric $\dY$.
Then 
\begin{enumerate}
\item $(\cGX, \fX, \dX)$ is \emph{Morita equivalent} to $(\cGY, \fY, \dY)$ if there exist clopen subsets $V \subset \fX$ and $W \subset \fY$, and a   homeomorphism $h \colon V \to W$ which conjugates $\cGX | V$ to $\cGY | W$.
\item $(\cGX, \fX, \dX)$ is \emph{Lipschitz equivalent} to $(\cGY, \fY, \dY)$ if the conjugation $h$ is Lipschitz.
\end{enumerate}
\end{defn}
  Morita equivalence is sometimes    called \emph{return equivalence} in the literature \cite{AO1995,Fokkink1991,CHL2013c}. 
  
  Morita equivalence is a basic notion   for the study of $C^*$-algebra invariants for foliation groupoids, as discussed by Renault \cite{Renault1980} and Connes \cite{Connes1994}. Lipschitz equivalence is a basic  notion   for the study of \emph{metric non-commutative geometry} \cite{Connes1994}. 

The strongest results for classification, up to Morita equivalence, have been obtained for   $1$-dimensional minimal matchbox manifolds. 
Fokkink  showed in his thesis \cite{Fokkink1991}  (see also Barge and Williams \cite{BargeWilliams2000})  that if $f_1, f_2$ are $C^1$-actions on $\mS^1$, each of which    have  Cantor minimal sets, then the induced minimal Cantor actions are Morita equivalent   if and only if  they have rotation numbers which are conjugate under the linear fractional action of $SL(2,\mZ)$ on $\mR$.  This implies   there are uncountably many non-homeomorphic minimal matchbox manifolds which embed as minimal sets for $C^1$-foliations of $\mT^2$. 
There is a higher-dimensional version of this result for torus-like matchbox manifolds, proved in \cite{CHL2013c}. See the papers \cite{BargeDiamond2001,BargeSwanson2007,BargeMartensen2011} for the classification of $1$-dimensional minimal matchbox manifolds embedded in   compact surfaces, which are necessarily not solenoids.  
 
In general, the classification problem  modulo orbit equivalence is unsolvable  for  the  pseudogroups associated to minimal matchbox manifolds of dimension $n \geq 2$, as  already for the normal solenoids with base manifold $\mT^n$ where $n \geq 2$, they are not classifiable. See \cite{Hjorth2000, KechrisMiller2004, Thomas2001, Thomas2003} for discussions of the undecidability of the Borel classification problem up to orbit equivalence.

 The advantage of considering Lipschitz equivalence of pseudogroup actions, is that while the equivalence is more refined, it can also be more practical to determine when two actions are not Lipschitz equivalent. We discuss   the difference between Morita and Lipschitz classification in the case of the weak solenoids, where there are a well-known criteria for Morita equivalence.

 First, we recall the criteria for when two weak solenoids are homeomorphic,  as described    in  \cite[Section~9]{ClarkHurder2013}, based on a result of  using a result of Mioduszewski \cite{Mioduszewski1963}.
 Assume that we are given two presentations, where   all spaces $\{M_{\ell} \mid \ell \geq 0\}$ and $\{N_{\ell} \mid \ell \geq 0\}$ are  compact oriented manifolds, and   all bonding maps are orientation-preserving coverings, 
\begin{equation}
\cP   =    \{p_{\ell+1} \colon M_{\ell+1} \to M_{\ell} \mid \ell \geq 0\} \quad, \quad 
\cQ   =    \{q_{\ell+1} \colon N_{\ell+1} \to N_{\ell} \mid \ell \geq 0\}
\end{equation}
which define   weak solenoids  $\cS_{\cP}$ and $\cS_{\cQ}$ as in \eqref{eq-presentationinvlim}, respectively. Choose basepoints $\ovx \in \cS_{\cP}$ and $\ovy \in \cS_{\cQ}$. We consider the special case where $M_0 = N_0$, as the more general case easily reduces to this one, and the key issues are more evident in this special case. 
Let $\Pi_{\ell}^{\cP} \colon \cS_{\cP} \to M_{\ell}$ denote the fibration map onto the factor $M_{\ell}$ for $\cS_{\cP}$, and $\Pi_{\ell}^{\cQ} \colon \cS_{\cQ} \to N_{\ell}$ that for $\cS_{\cQ}$.

We can assume that $x_0 = y_0$ in $M_0$, where $x_0 = \Pi_0^{\cP}(\ovx)$ and $y_0 = \Pi_0^{\cQ}(\ovy)$, then set $\cH_0 = \pi_1(M_0 , x_0)$, where we suppress the dependence on basepoints.  Define the   subgroups   $\cH_{\ell} \subset \cH_0$ which are the images of the groups $\pi_1(M_{\ell}, x_{\ell})$ under the maps $\ds (q_{\ell} )_{\#} $ associated to $\cP$, and let $\cG_{\ell} \subset \cH_0$ be the corresponding images of the groups $\pi_1(N_{\ell}, y_{\ell})$. Then we obtain two nested sequences of subgroups

 \begin{table}[htdp]
\begin{center}
\begin{tabular}{cccccccccc}
 $\subset$ & $\cH_{\ell+1}$ & $\subset$ &  $\cH_{\ell}$ &  $\subset$   & $\cdots$ & $\subset$ &     $\cH_{1}$ & $\subset$ &$\cH_{0}$ \\
 &  &  &   &    &   &   &     &  & $\parallel$ \\
 $\subset$ & $\cG_{\ell+1}$ & $\subset$ &  $\cG_{\ell}$ &  $\subset$   & $\cdots$ & $\subset$ &     $\cG_{1}$ & $\subset$ & $\cG_{0}$
\end{tabular}
\end{center}
\label{default}
\end{table}%
The proof of the following result can be found in the papers  \cite{McCord1965,Mioduszewski1963,Rogers1970,Schori1966}.
 \begin{thm}\label{thm-classifying1}
The  weak solenoids  $\cS_{\cP}$ and $\cS_{\cQ}$ are basepoint homeomorphic if and only if  there exists $\ell_0 \geq 0$ and $\nu_0 \geq 0$, such that  for every $\ell \geq \ell_0$ there exists $\nu_{\ell} \geq \nu_0$ with $\cG_{\nu_{\ell}} \subset \cH_{\ell}$, and
for every $\nu \geq \nu_0$ there exists $\ell_{\nu} \geq \ell_0$ with $\cH_{\ell_{\nu}} \subset \cG_{\nu}$. 
 \end{thm}
   The  condition on bonding maps   in Theorem~\ref{thm-classifying1} is called \emph{tower equivalence} of the subgroup chains.   
   
Let $\fX$ denote the fiber of  $\Pi_0^{\cP}$ over $\ovx$, and $\fY$ the fiber of $\Pi_0^{\cQ}$ over $\ovy$. Then the monodromy of the fibration $\Pi_0^{\cP}$ defines the actions of $\cH_0$ on $\fX$, and   the action of $\cH_0 = \cG_0$ on $\fY$ is defined by the monodromy of $\Pi_0^{\cQ}$.  Then results of   Clark, Lukina and the author  yield:

\begin{thm}[\cite{ClarkHurder2013}] \label{thm-classifying2}
If the weak solenoids  $\cS_{\cP}$ and $\cS_{\cQ}$ are basepoint homeomorphic, with   $M_0 = N_0$, then the holonomy actions of $\cH_0$ on $\fX$ and on $\fY$ are Morita equivalent.
\end{thm}

\begin{thm}[\cite{CHL2013c}] \label{thm-classifying3}
If the weak solenoids  $\cS_{\cP}$ and $\cS_{\cQ}$ have   base manifold $M_0 = N_0 = \mT^n$, and the holonomy actions of $\cH_0$ on $\fX$ and on $\fY$ are Morita equivalent, then  $\cS_{\cP}$ and $\cS_{\cQ}$ are basepoint homeomorphic.
\end{thm}
  
It follows that the classification problem for matchbox manifolds which are homeomorphic to a normal solenoid with base $\mT^n$,    reduces to the study of the Morita equivalence class of their holonomy pseudogroups, which by Theorem~\ref{thm-classifying1} reduces to a problem concerning the tower equivalence of subgroup chains in $\mZ^n$. The classification problem for subgroup chains    is not Borel,  for $n \geq 2$. 

In the case of classical Vietoris solenoids, where $M_0 = \mS^1$ and $\cH_0 = \mZ$,  the classification is much more straightforward. For each $\ell > 0$ there exists integers $m_{\ell} > 1$ and $n_{\ell} >  1$, defined recursively,  so that $\cH_{\ell} = \langle m_1 m_2 \cdots m_{\ell}\rangle \subset \mZ$, and 
 $\cG_{\ell} = \langle n_1 n_2 \cdots n_{\ell}\rangle \subset \mZ$. Let $P$ be the set of all prime factors of the integers $\{m_{\ell} \mid \ell > 0\}$, included with multiplicity, and let $Q$ be the same for the integers $\{n_{\ell} \mid \ell > 0\}$. For example, for the dyadic solenoid, the set $P = \{2,2,2,\ldots\}$ is   an infinite collection of copies of the prime $2$.  These infinite sets of primes $P$ and $Q$ are  ordered  by the sequence in which  they appear in the factorizations of the covering degrees $m_{\ell}$ and $n_{\ell}$.

 If the two sets $P$ and $Q$ are in \emph{bijective} correspondence, then it is an exercise to show that the  tower equivalence condition of Theorem~\ref{thm-classifying1} is satisfied for the presentations $\cP$ and $\cQ$, which yields the classification of  Vietoris solenoids up to homeomorphism by Bing     \cite{Bing1960} and McCord \cite{McCord1965} (see also Aarts and Fokkink \cite{AF1991}), and also the classification up to Morita equivalence of the associated minimal $\mZ$-actions on the Cantor set fibers. 
 
 However, for the metrics on the Cantor sections $\fX \subset \fM = \cS_{\cP}$ and $\fY \subset \fN = \cS_{\cP}$ as defined by the formula in \eqref{eq-canonicalmetric}, it is evident that if the bijection  $\sigma \colon P \leftrightarrow  Q$ permutes the elements by increasingly large degrees with respect to their ordering, then the induced map between the fibers,  $h_{\sigma} \colon \fX \cong \fY$, will not be Lipschitz.  This motivates introducing the following invariant of a tower of equivalences.

 Let $\cP$ and $\cQ$ be presentations with common base manifold $M_0$, and suppose there exists a tower equivalence between them. That is, there exists $\ell_0 \geq 0$ and $\nu_0 \geq 0$, such that  for every $\ell \geq \ell_0$ there exists $\nu_{\ell} \geq \nu_0$ with $\cG_{\nu_{\ell}} \subset \cH_{\ell}$, and
for every $\nu \geq \nu_0$ there exists $\ell_{\nu} \geq \ell_0$ with $\cH_{\ell_{\nu}} \subset \cG_{\nu}$. Define the \emph{displacement} of these indexing functions 
$\ell \mapsto \nu_{\ell}$ and $\nu \mapsto \ell_{\nu}$ to be
\begin{equation}
{\rm Disp}(\ell_{\nu}, \nu_{\ell}) = \max \left\{ \sup \left\{ |\ell_{\nu} - \nu| \ \mid \nu \geq \nu_0 \right\}~ , ~ \sup \left\{ |\nu_{\ell} - \ell| \ \mid \ell \geq \ell_0 \right\} \right\}
\end{equation}
If ${\rm Disp}(\ell_{\nu}, \nu_{\ell})  < \infty$, then we say that $\cP$ and $\cQ$ are \emph{bounded tower equivalent}. 
\begin{thm}\label{thm-Lipequivalent}
 Let $\cP$ and $\cQ$ be presentations with common base manifold $M_0$, and suppose there exists a tower equivalence between them, defined by maps $\ell \mapsto \nu_{\ell}$ and $\nu \mapsto \ell_{\nu}$. Let the fiber metrics be defined by the formula \eqref{eq-canonicalmetric} with $a_{\ell} = 3^{-\ell}$. Then the action of $\cH_0$ on the fiber $\fX$ of $\Pi_0^{\cP}$ is Lipschitz equivalent to the action of $\cH_0$ on the fiber $\fY$ of $\Pi_0^{\cQ}$ if and only if $\cP$ and $\cQ$ are bounded tower equivalent.
\end{thm}
The proof that ${\rm Disp}(\ell_{\nu}, \nu_{\ell}) < \infty$ implies Lipschitz equivalence for the metrics defined by \eqref{eq-canonicalmetric} with $a_{\ell} = 3^{-\ell}$ is  an exercise in the definitions, using the expression \eqref{eq-Galoisfiber} for the metric on the fibers. The converse direction, that Lipschitz equivalence implies bounded tower equivalence, follows from the works of Miyata and Watanabe \cite{MiyataWatanabe2002,MiyataWatanabe2003a}.

We give a simple example of Theorem~\ref{thm-Lipequivalent}, in the case of Vietoris solenoids. With the notation as above, suppose the the covering degrees $m_{\ell}$ for the presentation $\cP$ with base $M_0 = \mS^1$ are given by $m_{\ell} = 2$ for $\ell$ odd, and $m_{\ell} = 3$ for $\ell$ even. Let the covering degrees for the presentation $\cQ$ be given by the sequence $\{n_1, n_2, n_3, \ldots\} = \{2,3,2,2,3,2,2,2,2,3, \ldots\}$. In general, the $\ell$-th cover of degree $3$ is followed by $2^{\ell}$ covers of degree $2$. Then these two sequences are clearly tower equivalent, but their displacement is infinite. It follows that the matchbox manifolds  $\fM = \cS_{\cP}$ and $\fN = \cS_{\cP}$ are homeomorphic, but are not Lipschitz equivalent.


\end{document}